\documentclass[reqno,10pt]{amsart}
\usepackage{graphicx} % Required for inserting images
\usepackage{amsmath,amssymb,amsthm,mathtools}
\usepackage[a4paper,margin=3cm]{geometry}
\usepackage{amsfonts}
\usepackage{enumerate}
\usepackage{hyperref}
\usepackage{bbm}
\usepackage[inline]{enumitem}
\numberwithin{equation}{section}
\usepackage[foot]{amsaddr}
\usepackage{tikz}

\newcommand{\R}{\mathbb{R}}

\newcommand{\N}{\mathbb{N}}

\newcommand{\1}{\mathbbm{1}}

\title{An entropy for Boolean independence}

\author{Kewei Pan}

\begin{document}

\begin{abstract} 
    In this article, we introduce a notion of entropy for Boolean independence analogous to the Shannon's entropy and Voiculescu's free entropy through the study of asymptotic probabilities of the set of matrix approximates. To motivate this definition, we study two random matrix models exhibiting asymptotic Boolean independence, introduced respectively by Lenczewski and by C\'ebron--Gilliers. It turns out that the asymptotic probabilities of such matrix approximations can be characterized through large deviation principles for their empirical spectral measures. We prove that the associated rate functions coincide up to a scaling factor and attain their minimum at the Rademacher distribution, which plays the role of the Gaussian measure in Boolean probability. The logarithmic integral appearing in these rate functions is therefore identified as the Boolean entropy. We further show that this entropy is maximized, though not uniquely, by the Rademacher distribution and is monotone along the Boolean central limit theorem.
\end{abstract}

\maketitle

\tableofcontents

\newtheorem{theorem}{Theorem}[section] 
\newtheorem{lemma}[theorem]{Lemma}
\newtheorem{corollary}[theorem]{Corollary}
\newtheorem{proposition}[theorem]{Proposition}
\theoremstyle{remark}
\newtheorem{remark}[theorem]{Remark}
\newtheorem{definition}[theorem]{Definition}
\newtheorem{example}[theorem]{Example}
\newtheorem{notation}[theorem]{Notation}

\section{Introduction}
\subsection{Background}
\noindent\textbf{Noncommutative independence.} A noncommutative probability space is a pair $(\mathcal{A},\varphi)$, where $\mathcal{A}$ is a unital algebra with unit $1_\mathcal{A}$ and $\varphi:\mathcal{A}\to\mathbb{C}$ is a linear functional such that $\varphi(1_\mathcal{A})=1$. The elements of $\mathcal A$ are regarded as noncommutative random variables, and the distribution of $a\in\mathcal A$ is encoded by its moments $\{\varphi(a^n)\}_{n\ge 1}$. When $\mathcal A$ is a $C^*$-algebra and $a=a^*$, the distribution naturally extends to a compactly supported probability measure $\mu_a\in\mathcal M(\mathbb R)$, called the distribution of $a$. The notion of a noncommutative probability space generalizes that of a classical probability space. It is natural to ask whether there exist genuinely noncommutative notions of independence other than the classical one. A rich and nontrivial noncommutative probability theory was developed in the mid-1980s by Voiculescu \cite{Voiculescu1985symmetries}, who introduced the notion of \emph{freeness} in his study of free group factors. Subsequently, Speicher \cite{speicher1996universal} introduced another notion of independence called \emph{Boolean independence} and further showed that under suitable assumptions only three universal notions of independence arise in noncommutative probability spaces: \emph{tensor (classical), free, and Boolean independence}. The definitions are as follows. Let $(\mathcal{A},\varphi)$ be a noncommutative probability space. Given sub-algebras $(\mathcal{A}_k)_{k\geq 1}\subset \mathcal{A}$ and $a_j\in \mathcal{A}_{i(j)}\,(j=1,\ldots,n),$ we define the associated partition $\pi=\{V_1,\ldots,V_p\}\in\mathrm{P}(n)$ by the relation: $j\sim_{\pi}l$ if $i(j)=i(l)$. For any partition $\sigma=\{W_1,\ldots,W_q\}\in \mathrm{P}(n),$ we put
    \begin{equation*}
        \varphi_{\sigma}(a_1\cdots a_n):=\varphi(\prod_{j\in W_1}a_j)\cdots \varphi(\prod_{j\in W_q}a_j).
    \end{equation*}

    \begin{enumerate}[leftmargin=20pt,label=(\arabic*)]
        \item We say $(\mathcal{A}_k)_{k\in\mathbb{N}}$ are \textbf{tensor} independent if for any $\{a_1,\ldots,a_n\}\in\bigcup_{k\in\mathbb{N}}\mathcal{A}_k$ (denote $\pi\in\mathrm{P}(n)$ the associated partition), 
        \begin{equation*}
            \varphi(a_1\cdots a_n)= \varphi_{\pi}(a_1\cdots a_n).
        \end{equation*}
        This notion can be viewed as a noncommutative extension of classical independence. We shall also refer to it as classical independence when we restrict to the commutative setting.
        \item We say $(\mathcal{A}_k)_{k\in\mathbb{N}}$ are \textbf{free} if for any $\{a_1,\ldots,a_n\}\in\bigcup_{k\in\mathbb{N}}\mathcal{A}_k$
        such that $i(j)\neq i(j+1),\ 1\leq j\leq n-1,$ 
        \begin{equation*}
            \forall j,\ \varphi(a_j)=0 \Longrightarrow\varphi(a_1,\ldots,a_n)=0.
        \end{equation*}
        \item We say $(\mathcal{A}_k)_{k\in\mathbb{N}}$ are \textbf{Boolean} independent if for any $\{a_1,\ldots,a_n\}\in\bigcup_{k\in\mathbb{N}}\mathcal{A}_k$ such that $i(j)\neq i(j+1),\ 1\leq j\leq n-1$ and $a_i\ne 1_\mathcal{A}$,
        \begin{equation*}
            \varphi(a_1\cdots a_n)=\prod_{j=1}^n\varphi(a_j).
        \end{equation*}
    \end{enumerate}
We say that two random variables $a,b\in\mathcal{A}$ are independent with respect to a given notion of independence if the subalgebras $\operatorname{Alg}(a)$ and $\operatorname{Alg}(b)$ generated by $a$ and $b$ are independent in that sense.

It turns out that each notion of independence admits an associated family of cumulants, and the corresponding independence is characterized by the vanishing of mixed cumulants. Moreover, these cumulants can be recovered from moments through Möbius inversion on different lattices of partitions. The classical cumulants are known to be indexed by all partitions (P). In contrast, freeness is governed by non-crossing partitions (NC), see, e.g., \cite{speicher1994multiplicative}. In the Boolean setting \cite{speicher1996boolean}, the corresponding combinatorial structure is given by interval partitions (I).

One can also study the central limit theorem in each of these settings. More precisely, consider a sequence of identically distributed self-adjoint random variables (i.e., having the same moments) $\{a_i\}_{i\geq1}\subset(\mathcal{A},\varphi)$ such that $a_1,a_2,\ldots$ are independent with respect to a given notion of independence. Furthermore, assume that $\varphi(a_1)=0,\,\varphi(a_1^2)=1$ and define
\[
s_n=\frac{a_1+\cdots+a_n}{\sqrt{n}}.
\]
The goal is to determine the limiting probability measure $\mu_{\mathrm{clt}}\in\mathcal{M}(\mathbb{R})$ such that
\[
\lim_{n\to\infty}\varphi(s_n^k)=m_k(\mu_{\mathrm{clt}}),
\qquad k\geq1,
\]
where $m_k(\mu_{\mathrm{clt}})$ denotes the $k$-th moment of $\mu_{\mathrm{clt}}$. Using cumulants, one can compute the limiting moments explicitly, which are equal to the enumeration of pairings within the corresponding classes of partitions.

\begin{table}[ht]
      \label{clt}
      \centering
      \begin{tabular}{c|c|c|c}
      Independence & tensor   &  free &	 Boolean	\\
      \hline
      Lattice & P & NC & I\\
      \hline
      $m_k(\mu_{\mathrm{clt}})$ & $\#\mathrm{P}_2(k)$ & $\#\mathrm{NC}_2(k)$ & $\#\mathrm{I}_2(k)$\\
      \hline
      $\mu_{\mathrm{clt}}$ & $\gamma_1\sim\mathcal{N}(0,1)$ &  $\sigma(dx)\sim \sqrt{4-x^2}\,dx$ & $\mathrm{b}=\frac{1}{2}\delta_{-1}+\frac{1}{2}\delta_1$
      \end{tabular}
      \caption{Parallel structures among three universal notions of independence. The free CLT was first proved by Voiculescu \cite{voiculescu1986addition} using an analytical approach, while a combinatorial proof via free cumulants was later given by Speicher \cite{speicher1990new}; The Boolean CLT was established by Speicher and Woroudi \cite{speicher1996boolean}.}
      \label{table:parallel independence}
\end{table}

\medskip
\noindent\textbf{Classical entropy.}
Let $X=(X_1,\ldots,X_d)$ be a classical random vector with distribution
$\mu\in\mathcal{M}(\mathbb{R}^d)$ that is compactly supported. The Shannon's (differential) entropy of $X$ is defined by
\begin{equation*}
    S(X)=S(\mu)
    :=
    -\int_{\mathbb{R}^d} f(x)\log f(x)\,dx,
    \qquad
    \text{if }
    \frac{d\mu}{dx}=f,
\end{equation*}
and $S(\mu)=-\infty$ otherwise.

The notion of entropy originated in the work of Boltzmann in statistical mechanics to quantify the asymptotic behavior of large systems. From a probabilistic perspective, entropy can be recovered from the following asymptotic probability as $n\to\infty$, $r\to\infty$ and $\varepsilon\to0$,
\begin{equation}\label{eq:classical asym prob}
\frac{1}{n}\log
\gamma_d^n\!\left( (x_i)_{i=1}^d\in (\mathbb{R}^n)^d:
\left|
\mathrm{tr}_n(D_{i_1}\cdots D_{i_k})
-
\mathbb{E}[X_{i_1}\cdots X_{i_k}]
\right|
\le \varepsilon,\;
k\le r
\right),
\end{equation}
where $D_i=\operatorname{diag}\bigl(x_i(1),\ldots,x_i(n)\bigr)$ for $1\le i\le d,$ and $\gamma_d^n$ denotes the $n$-fold product of the standard Gaussian measure on $\mathbb{R}^d$. Equivalently, write $x(j)=(x_1(j),\ldots,x_d(j))$,
we associate the empirical measure
\[
\frac{1}{n}\sum_{j=1}^n \delta_{x(j)}\in \mathcal{M}(\mathbb{R}^d).
\]
In this formulation, \eqref{eq:classical asym prob} measures the proximity of the empirical measure to $\mu$ by requiring that the mixed moments up to order $r$ differ from those of $\mu$ by at most $\varepsilon$. Indeed, this is an application of the Sanov theorem \cite{sanov1961probability}, which says the empirical measure defined above satisfies a large deviation principle with speed $n$ and a rate function given by
\[
H(\mu\mid\gamma_d) =
\int_{\mathbb{R}^d}
\frac{|x|^2}{2}\,d\mu(x)
-
S(\mu)
+\mathrm{const},
\]
which is just the limit of \eqref{eq:classical asym prob}. This rate function is also known as the relative entropy and admits a unique minimizer, given by $\gamma_d$. Consequently, the Gaussian measure maximizes entropy among all probability measures with a prescribed second moment. Remarkably, when \(d=1\), along the classical CLT, not only is the entropy bounded above by its maximal value $S(\gamma_1)$, but also it satisfies the monotonicity property
\[
S(s_n)\nearrow S(\gamma_1).
\]
The convergence to $S(\gamma_1)$ was proved by Barron \cite{andrew1986entropy} and the problem of monotonicity of entropy dates back to Shannon (see \cite{shannon1948mathematical,stam1959some,lieb1978proof}), and remains open for many years and was eventually resolved by the authors in \cite{artstein2004solution}.

\medskip
\noindent\textbf{Free entropy.} 
Given a random variable $a\in\mathcal A$ with distribution $\mu\in\mathcal{M}(\mathbb{R})$ that is compactly supported, Voiculescu \cite{voiculescu1993analogues} observed that, for a sequence of unitarily invariant self-adjoint random matrices $A_N$ whose empirical spectral measures converge to $\mu$, the Boltzmann-Gibbs entropy of $A_N$ formally behaves, in the large-$N$ limit, like
\[
\Sigma(a)=\Sigma(\mu):=
\iint \log |x-y|\,d\mu(x)d\mu(y).
\]
He therefore introduced $\Sigma$ as the \emph{free entropy}, which coincides with minus the logarithmic energy of $\mu$. Nevertheless, this one-dimensional quantity alone does not explain the adjective ``free.'' Motivated by the asymptotic freeness of large independent random matrices established in \cite{voiculescu1991limit}, Voiculescu  \cite{voiculescu1994analogues} subsequently introduced a microstates free entropy, which extends the above definition to multivariate case. One of its fundamental properties is its additivity for free random variables, mirroring the additivity of classical entropy for independent random variables. For a $d$-tuple $a=(a_1,\ldots,a_d)\subset (\mathcal A,\varphi),$
the basic idea is to approximate $a$ by a $d$-tuple of independent random matrices and study the asymptotic probability
\begin{equation}\label{eq:free asym prob}
\frac{1}{N^2}\log
\mathbb{P}\!\left(
\{A_i\}_{i=1}^d\in M_N(\mathbb{C})^d:
\left|
\operatorname{tr}_N(A_{i_1}\cdots A_{i_k})
-
\varphi(a_{i_1}\cdots a_{i_k})
\right|
\le \varepsilon,\;
k\le r
\right).
\end{equation}
This construction is parallel to \eqref{eq:classical asym prob}, with diagonal matrices replaced by Hermitian matrices. In general, no explicit formula is known for the multivariate free entropy. However, in the one-dimensional case, the free entropy can also be computed through the large deviation theory of random matrices.

Consider the Gaussian Unitary Ensemble (GUE),
\[
G_N\sim \exp\!\left(-\frac{1}{2}\operatorname{Tr}(G_N^2)\right)dG_N,
\]
where $dG_N$ denotes the Lebesgue measure on the space of $N\times N$ Hermitian matrices. Let $X_N=N^{-1/2}G_N$ and denote empirical spectral measure of $X_N$
\[
\mu_{X_N} =
\frac1N\sum_{i=1}^N\delta_{\lambda_i}.
\]

Now let $a$ be a random variable with compactly supported distribution $\mu\in\mathcal{M}(\mathbb{R})$. Specializing \eqref{eq:free asym prob} to the one-dimensional setting yields
\[
\frac{1}{N^2}
\log
\mathbb{P}
\Bigl(X_N:
|m_k(\mu_{X_N})-m_k(\mu)|
\le \varepsilon,\;
k\le r
\Bigr),
\]
which can be reformulated in terms of the large deviation principle for the empirical spectral measure $\mu_{X_N}$. Using the explicit joint density of the eigenvalues,
\begin{equation}\label{eq:GUE density}
dP_N(\lambda_1,\ldots,\lambda_N)
=
\frac1{Z_N}
\prod_{i<j}|\lambda_i-\lambda_j|^2
\exp\!\left(
-N\sum_{i=1}^N\frac{\lambda_i^2}{2}
\right)
\prod_{i=1}^N d\lambda_i,
\end{equation}
Ben Arous and Guionnet \cite{arous1997large} established a large deviation principle for $\mu_{X_N}$ with speed $N^2$ and rate function
\begin{equation}\label{eq:rate function GUE}
   \int \frac{x^2}{2}\,d\mu(x)
-\Sigma(\mu)
-\frac34.
\end{equation}
The unique minimizer of the rate function is the semicircle law $\sigma$, which plays the same role in free probability as the Gaussian distribution does in the classical setting.

It turns out that free entropy shares many analogous properties with classical entropy. For instance, the semicircle law maximizes free entropy under a prescribed variance constraint, and $\Sigma(s_n)$ is non-decreasing along the free CLT \cite{shlyakhtenko2007}. The convergence of $\Sigma(s_n)$ to the semicircular free entropy $\Sigma(\sigma)$ was later proved by Wang \cite{Wang2010locallimit}.

\medskip
The goal of this work is to develop a Boolean analogue of entropy within the paradigm established by classical and free entropy. In both settings, the underlying notion of independence is encoded by suitable random matrix models: diagonal matrices for classical independence and Hermitian matrices for freeness. Since these matrices become asymptotically independent in their respective senses, the associated asymptotic probabilities naturally lead to the notions of classical and free entropy. Motivated by this principle, we introduce two random matrix models exhibiting asymptotic Boolean independence and investigate the corresponding asymptotic probabilities. We only focus on one-dimensional case in this article, and propose the resulting quantity as a natural candidate for Boolean entropy.

\subsection{Models}
It turns out that Boolean independence can also arise from large random matrix models, but typically in a degenerate regime where most eigenvalues collapse to zero. 

\medskip
\noindent\textbf{Symmetric Block model.} The first model was introduced by Lenczewski \cite{lenczewski2014limit}. He considered two independent GUE matrices $\{G_n(1),G_n(2)\}$ (more general models are studied in the original paper) and their associated symmetric block matrices
\begin{equation*}
    T_{p,n}(i)=
    \begin{pmatrix}
        O & G_{p,q}(i)\\
        G_{p,q}(i)^* & O
    \end{pmatrix},\qquad i=1,2
\end{equation*}
where $p+q=n$ and $G_{p,q}(i)\in M_{p\times q}(\mathbb{C})$ has i.i.d. complex Gaussian entries with variance $1$. Equivalently, the law of $T_{p,n}$ can be written as
\begin{equation*}
    T_{p,n}\sim
    \exp\left(
        -\frac12 \mathrm{Tr}_p(T^2)
    \right)\,dT,
\end{equation*}
where $\mathrm{Tr}_p$ denotes the partial trace on the first block,
\begin{equation*}
    \mathrm{Tr}_p(M)
    =
    \sum_{i=1}^{p}
    \langle Me_i,e_i\rangle,
\end{equation*}
with $\{e_1,\ldots,e_n\}$ the canonical basis of $\mathbb{C}^n$, and $dT$ denotes the Lebesgue measure on the space of Hermitian matrices of the form
\begin{equation*}
    T=
    \begin{pmatrix}
        O & A\\
        A^* & O
    \end{pmatrix},
    \qquad
    A\in M_{p\times q}(\mathbb{C}).
\end{equation*}
Assume that $p=p_n$ depends on $n$ in such a way that $p_n\to\infty$ and $p_n/n\to0$ as $n\to\infty$. Lenczewski showed that the family
\[
\left\{\frac{1}{\sqrt n}T_{p,n}(1),\,\frac{1}{\sqrt n}T_{p,n}(2)\right\}
\]
is asymptotically Boolean independent with respect to the normalized partial trace $\mathrm{tr}_p=p^{-1}\mathrm{Tr}_p$. We shall refer to this random matrix ensemble as the Symmetric Block Model. 

Motivated by \eqref{eq:free asym prob}, we define the approximation set by replacing the normalized trace $\operatorname{tr}_N$ with the normalized partial trace $\operatorname{tr}_p$, and the independent Hermitian matrices $A_i$ with independent symmetric block matrices $n^{-1/2}T_{p,n}(i)$. In particular, in the one-dimensional setting, the central object of study becomes the moments of $T_{p,n}$ with respect to $\operatorname{tr}_p$.

To this end, let $\{s_i\}_{i=1}^p$ denote the singular values of $n^{-1/2}G_{p,q}$. Note that the eigenvalues of $n^{-1/2}T_{p,n}$ come in symmetric pairs,
\begin{equation*}
    \mathrm{Spec}(n^{-1/2}T_{p,n})
    =
    \left\{
    \pm s_1,\ldots,\pm s_p,
    0,\ldots,0
    \right\}.
\end{equation*}
Moreover, since $W(p,n)=n^{-1}G_{p,q}G_{p,q}^*$
is the $p\times p$ complex Wishart matrix, its eigenvalue density is given (up to a normalization constant) by
\begin{equation}\label{eq:wishart eigen density}
    \prod_{1\leq i<j\leq p}|\lambda_i-\lambda_j|^2
    \cdot
    \prod_{i=1}^p\lambda_i^{n-p}
    \cdot
    \prod_{i=1}^p e^{-n\lambda_i}
    \,d\lambda_i,
\end{equation}
with respect to the Lebesgue measure on $\mathbb{R}_+^p$ (see, e.g., \cite[Section 5.5]{hiai2000semicircle}). By the change of variables $\lambda_i=s_i^2$, we obtain the density of $(s_1,\ldots,s_p)\in\mathbb{R}_+^p$:
\begin{equation}\label{eq:singular density}
    dW_{p,n}(s_1,\ldots,s_p)
    =
    \frac{1}{D_{p,n}}
    \prod_{1\leq i<j\leq p}|s_i^2-s_j^2|^2
    \cdot
    \prod_{i=1}^p s_i^{2(n-p)+1}
    \cdot
    \prod_{i=1}^p e^{-ns_i^2}
    \,ds_i.
\end{equation}
Finally, a direct computation shows that
\begin{equation}\label{symmetric block equality}
    \mathrm{tr}_p((n^{-1/2}T_{p,n})^k)
    =
    \frac{1}{2p}\mathrm{Tr}((n^{-1/2}T_{p,n})^k).
\end{equation}
Therefore, the asymptotic behavior of the model is completely encoded by the symmetric averaged singular value distribution
\begin{equation}\label{def: sym block esd}
    \hat{\mu}_p^T
    :=
    \frac{1}{p}
    \sum_{i=1}^{p}
    \left(
    \frac{1}{2}\delta_{s_i}
    +
    \frac{1}{2}\delta_{-s_i}
    \right)
    \in
    \mathcal{M}^{\mathrm{sym}}(\mathbb{R}),
\end{equation}
where $\mathcal{M}^{\mathrm{sym}}(\mathbb{R})$ denotes the space of symmetric probability measures on $\mathbb{R}$, namely those satisfying $\mu(B)=\mu(-B)$ for every Borel set $B$. 

\medskip
\noindent\textbf{Conditioned GUE model.} Recently, C{\'e}bron and Gilliers \cite{cebron2024asymptotic} introduced the \emph{Vortex model}. Let $\{v_N\}_{N\geq1}$ be a sequence of deterministic unit vectors in $\mathbb{C}^N$, and define a linear functional $\psi^{v_N}:M_N(\mathbb{C})\to\mathbb{C}$ by
\begin{equation*}
    \psi^{v_N}(M)
    =
    \langle Mv_N,v_N\rangle,
    \qquad
    M\in M_N(\mathbb{C}).
\end{equation*}
Furthermore, let $\{\tilde{U}_N\}_{N\geq1}$ be a sequence of Haar-distributed $N\times N$ unitary matrices conditioned to leave $v_N$ invariant, that is, $\tilde{U}_Nv_N=v_N$. They proved that Boolean independence emerges asymptotically between $A_N$ and the randomly rotated matrix $\tilde{U}_N^*B_N\tilde{U}_N$ with respect to the state $\psi^{v_N}$, provided that $\{A_N,B_N\}_{N\geq1}$ are sequences of deterministic self-adjoint matrices with uniformly bounded operator norms and whose spectra asymptotically collapse to $0$. More precisely, the empirical spectral measures of $A_N$ and $B_N$ converge to $\delta_0$ under the normalized trace $\mathrm{tr}_N$ as $N\to\infty$.

Recall that for the normalized GUE matrix $X_N=N^{-1/2}G_N$, one has the spectral decomposition
\begin{equation*}
    X_N=U_N^*\mathrm{diag}(\lambda_1,\ldots,\lambda_N)U_N.
\end{equation*}
where the unitary matrix $U_N$ is Haar distributed and independent of the eigenvalues $\{\lambda_i\}_{i=1}^N$. Moreover, the joint density of the eigenvalues of $X_N$ is given by \eqref{eq:GUE density}. Let $1\leq M\leq N$, Conditioned on $\lambda_{M+1}=\cdots=\lambda_N=0$, the joint density of $(\lambda_1,\ldots,\lambda_M)$ is (up to a normalization constant),
\begin{equation}\label{eq:con GUE density}
    \prod_{i=1}^M\lambda_i^{2(N-M)}
    \prod_{1\leq i<j\leq M}|\lambda_i-\lambda_j|^2
    \prod_{i=1}^M e^{-N\frac{\lambda_i^2}{2}}d\lambda_i.
\end{equation}
We choose $v_N=M^{-1/2}(e_1+\cdots+e_M),$
where $\{e_i\}_{i=1}^M$ are the canonical basis vectors of $\mathbb{R}^N$, and further condition $U_N$ to leave $v_N$ invariant. Since $U_N$ is independent of the eigenvalues, the resulting Conditioned GUE matrix has the decomposition
\begin{equation*}
    \tilde{X}_N
    =
    \tilde{U}_N^*
    \mathrm{diag}(\lambda_1,\ldots,\lambda_M,0,\ldots,0)
    \tilde{U}_N,
\end{equation*}
where $\tilde{U}_Nv_N=v_N$ and $(\lambda_1,\ldots,\lambda_M)$ follows the density in \eqref{eq:con GUE density}. Therefore, this Conditioned GUE model is a special case of the Vortex model introduced above. Furthermore, if we assume that $M=M_N$ depends on $N$ such that $M\to +\infty,\ M/N\to 0$ as $N\to +\infty$, two independent copies of $\tilde{X}_N$ are asymptotically Boolean independent.

By the choice of the vector state $\psi^{v_N}$, we have for any integer $k\ge 1$,
\begin{equation*}
    \psi^{v_N}(\tilde{X}_N^k)
    =
    \frac{1}{M}\sum_{i=1}^M\lambda_i^k
    =
    \int x^k dL_M(x),
\end{equation*}
where
\begin{equation}\label{def:con GUE esd}
    L_M:=\frac{1}{M}\sum_{i=1}^M\delta_{\lambda_i}.
\end{equation}
In the same spirit, studying the asymptotic probability of the Conditioned GUE model $\tilde{X}_N$ with respect to the state $\psi^{v_N}$ can be reduced to study the large deviation principle of the empirical spectral measure $L_M$.

\subsection{Main results}
In light of these results, the main goal of this article is to investigate the large deviation principles for the Symmetric Block model and the Conditioned GUE model. A priori, the rate functions arising from different models may differ. We prove that these rate functions indeed coincide, thus providing a natural candidate for the definition of Boolean entropy, which is independent of the choice of the underlying model.

First we recall the definition of large deviation principle: let $X$ to be a Polish space (that is a complete separable metric space). We recall that a function $f: X\to \mathbb{R}$ is lower semicontinuous if the level sets $\{x: f(x)\leq C\}$ are closed for any constant $C$.
\begin{definition}\label{LDP}
    A sequence $(\mu_N)_{N\in\mathbb{N}}$ of probability measures on $X$ satisfies a large deviation principle with speed $a_N$ (going to infinity with $N$) and rate function $I$ iff
\begin{enumerate}[leftmargin=20pt, label=(\alph*)]
    \item $I: X\to [0,\infty]$ is lower semicontinuous.
    \item For any open set $O\in X, \displaystyle\liminf_{N\to \infty}\frac{1}{a_N}\log\mu_N(O)\geq -\inf_O I.$
    \item For any closed set $F\in X, \displaystyle\limsup_{N\to \infty}\frac{1}{a_N}\log\mu_N(F)\leq -\inf_F I.$
\end{enumerate}
Moreover, a rate function is good if for any $M\geq 0,$ the set $\{x\in X: I(x)\leq M\}$ is compact.
\end{definition}

We endow $\mathcal{M}(\mathbb{R})$ with the weak topology, which is compatible with the Lipschitz bounded metric:
    \begin{equation}\label{Lipschitz metric}
        d_{BL}(\mu,\nu)=\sup_{f\in\mathcal{F}_{LU}}\left|\int fd\mu-\int fd\nu\right|,
    \end{equation}
where $\mathcal{F}_{LU}$ is the class of Lipschitz continuous functions $f: \mathbb{R} \to \mathbb{R}$ with Lipschitz constant at most $1$ and uniform bound $1$.

First, we present a large deviation principle for the Symmetric Block model.
\begin{theorem}\label{thm:LDP sym block}
    In the regime $p=p_n\to\infty$ and $p/n\to0$ as $n\to\infty$, under the law \eqref{eq:singular density}, the symmetric averaged singular value distribution $\hat{\mu}_p^T$ defined in \eqref{def: sym block esd} satisfies a large deviation principle on $\mathcal{M}^{\mathrm{sym}}(\mathbb{R})$ with speed $pn$ and good rate function
    \begin{equation}\label{eq:rate function sym block}
       I^{\mathrm{sym}}(\mu):=\int (x^2-\log x^2)\,d\mu(x)-1.
    \end{equation}
    Moreover, $I^{\mathrm{sym}}$ admits a unique minimizer, given by the Rademacher distribution $\mathrm{b}$.
\end{theorem}

We also establish a large deviation principle for the Conditioned GUE model.
\begin{theorem}\label{thm:LDP con GUE}
   In the regime $M=M_N\to\infty$ and $M/N\to0$ as $N\to\infty$, under the law \eqref{eq:con GUE density}, the empirical spectral measure $L_M$ defined in \eqref{def:con GUE esd} satisfies a large deviation principle on $\mathcal{M}(\mathbb{R})$ with speed $NM$ and good rate function
   \begin{equation}\label{eq:rate function con GUE}
      I(\mu):=\int \left(\frac{1}{2}x^2-\log x^2\right)d\mu(x)-\frac{1}{2}.
   \end{equation}
   The minimizer of $I$ is not unique. Indeed, $p\delta_{-\sqrt{2}}+(1-p)\delta_{\sqrt{2}}$ minimizes $I$ for every $p\in[0,1]$.
\end{theorem}
As a corollary, the empirical spectral measure \(L_M\) is asymptotically concentrated near the set
\[
\left\{p\delta_{-\sqrt{2}}+(1-p)\delta_{\sqrt{2}}: p\in[0,1]\right\}.
\]
By symmetry considerations, one expects \(L_M\) to converge, at least in expectation, to $\frac{1}{2}\delta_{-\sqrt{2}}+\frac{1}{2}\delta_{\sqrt{2}}.$
It is therefore natural to ask whether this convergence also holds almost surely. This motivates us to establish a refined large deviation principle for the rescaled empirical spectral measures
\[
\mu_{\alpha,M}
=
\frac{1}{M}\sum_{i:\lambda_i\ge 0}\delta_{\alpha_i},
\qquad
\mu_{\beta,M}
=
\frac{1}{M}\sum_{i:\lambda_i<0}\delta_{-\beta_i},
\]
where
\[
\alpha_i=\Lambda_M^{-1}(\lambda_i-\sqrt{2}),
\qquad
\beta_i=\Lambda_M^{-1}(\lambda_i+\sqrt{2}),
\]
and $\{\Lambda_M\}_M$ is a sequence of constants depending on $M$. With a suitable choice of $\Lambda_M$ and in an appropriate Polish space, we prove that the pair $(\mu_{\alpha,M},\mu_{\beta,M})$
satisfies a large deviation principle (see Theorem \ref{thm:LDP rescaled con GUE} for details). As a consequence, we obtain
\[
\frac{\#\{i:\lambda_i\ge0\}}{M}\longrightarrow\frac12
\]
almost surely, which together with Theorem \ref{thm:LDP con GUE} determines the limiting distribution of $L_M$.

In the light of \cite{arous1997large}, the same argument applies to a more general class of potentials. Consider the Gibbs measure on $N\times N$ Hermitian matrices
$$
dP_N^V(H)\sim \exp(-N\mathrm{Tr}(V(H)))\,dH,
$$
where $V:\mathbb{R}\to\mathbb{R}$ satisfies $V(x)-\log|x|\to+\infty$ as $|x|\to+\infty$. Then the law of $(\lambda_1,\ldots,\lambda_M)$ conditioned on $\lambda_{M+1}=\cdots=\lambda_N=0$ has a density proportional to
\begin{equation}\label{eq:con Gibbs potential density}
    \prod_{i=1}^M\lambda_i^{2(N-M)}
    \prod_{1\leq i<j\leq M}|\lambda_i-\lambda_j|^2
    \prod_{i=1}^M e^{-NV(\lambda_i)}\,d\lambda_i.
\end{equation}
Further extensions are discussed in Section 3.3, and here we present the corresponding large deviation principle. In the regime $M=M_N\to\infty$ and $M/N\to0$ as $N\to\infty$, under the law \eqref{eq:con Gibbs potential density}, the empirical spectral measure $L_M$ defined in \eqref{def:con GUE esd} satisfies a large deviation principle on $\mathcal{M}(\mathbb{R})$ with speed $NM$ and good rate function
\begin{equation}
    I^V(\mu):=\int \left(V(x)-\log x^2\right)d\mu(x)-C_V.
\end{equation}
There is no difficulty extending our results to the conditioned $\beta$-ensemble (see Corollary \ref{cor:LDP con beta ensemble}). One may also consider the Symmetric Block model under the Gibbs measure $P_N^V$ on matrices. However, such a generalization remains out of reach at present. The main obstacle is the lack of an explicit matrix ensemble realizing the corresponding Gibbs distribution.

Based on the previous large deviation results, and in analogy to the construction of classical entropy and Voiculescu's free entropy, we propose the logarithmic integral
\begin{equation}
    \Gamma(a)=\Gamma(\mu):=\int \log |x|\,d\mu(x),\qquad \text{if $a$ has law $\mu$}
\end{equation}
as a candidate for (microstates) Boolean entropy. In parallel with the classical and free settings, several aspects remain to be explored in order to further justify the notion of Boolean entropy, most notably the additivity property for Boolean independent random variables. These questions are left for future work. In this article, we restrict ourselves to studying the maximality and monotonicity of \(\Gamma\) along the Boolean central limit theorem.

\begin{theorem}\label{thm:max and mono}
    We have the following:
    \begin{enumerate}[leftmargin=20pt, label=(\alph*)]
        \item Among the set of probability measures with $m_2(\mu)=1$, $p\delta_{-1}+(1-p)\delta_1$ maximizes $\Gamma(\mu)$ for any $p\in[0,1]$.
        
        \item Let $\{a_i\}_{i\ge1}\subset(\mathcal{A},\varphi)$ be a sequence of identically distributed and Boolean independent random variables. If $\Gamma(a_1)>-\infty$, then
        \begin{equation*}
            \Gamma\left(\frac{a_1+\cdots+a_n}{\sqrt{n}}\right)
        \end{equation*}
        is non-decreasing in $n$.
    \end{enumerate}
\end{theorem}
\begin{remark}
    In particular, the Boolean entropy is non-decreasing along the Boolean CLT. However, it remains unclear whether $\Gamma$ converges to its maximal value $\Gamma(\mathrm{b})$ along the Boolean CLT. One may also define Boolean entropy on the subspace \(\mathcal{M}^{\mathrm{sym}}(\mathbb{R})\). A key advantage of this restriction is that, under a prescribed second moment, the Boolean entropy \(\Gamma\) admits a unique maximizer, given by the Rademacher distribution, which is the Boolean analogue of the Gaussian distribution (see Table \ref{table:parallel independence}). This parallels the maximality properties of both classical entropy and free entropy. We further remark that the notion $\Gamma$ is not the unique notion for Boolean entropy, at the level of justification on two properties in Theorem \ref{thm:max and mono}. Indeed, one can easily construct a real symmetric function $\omega$ such that $\frac{1}{2}x^2-\omega(x)$ has minimizers at $\pm 1$, and thus $\int \omega(x)\,d\mu(x)$ is maximized by the Rademacher distribution in $\mathcal{M}^{\mathrm{sym}}(\mathbb{R})$. Moreover, there exists other functional satisfying the same monotonicity (see Corollary \ref{cor:general h along boolean clt} and Remark \ref{rem:existence of other functional for monotonicity}).
\end{remark}

\noindent\textbf{Organization of the paper.}
The proof of Theorem \ref{thm:LDP sym block} is given in Section 2, where we also explain its connection with the large deviations of the Wishart model. Section 3 is devoted to the study of the Conditioned GUE model: The proof of Theorem \ref{thm:LDP con GUE} follows a similar argument to that of Theorem \ref{thm:LDP sym block}; we therefore only provide a sketch of the proof, the main focus is the refined large deviation principle (Theorem \ref{thm:LDP rescaled con GUE}) and its proof are established in Section 3.2, while several extensions are discussed in Section 3.3. Finally, we prove Theorem \ref{thm:max and mono} in Section 4.

\medskip
\noindent\textbf{Acknowledgements.} The author would like to thank Guillaume Cébron and Djalil Chafaï for suggesting this topic and for several helpful discussions. The author is also grateful to the referees for their time and valuable comments.

\section{LDP for the Symmetric Block model}
\subsection{Proving Theorem \ref{thm:LDP sym block}}
Recall that
\begin{equation*}
    \hat{\mu}_p^T:=\frac{1}{2p}\sum_{i=1}^{p}\left(\delta_{s_i}+\delta_{-s_i}\right).
\end{equation*}
Due to the symmetry, it suffices to prove the large deviation for $\hat{\mu}_{p,+}^T:=p^{-1}\sum_{i=1}^p\delta_{s_i}\in\mathcal{M}(\R_+),$ where $(s_1,\ldots,s_p)\in\mathbb{R}_+^p$ follows the density:
\begin{equation*}
    dW_{p,n}(s_1,\ldots,s_p)=\frac{1}{D_{p,n}}\prod_{1\leq i<j\leq p}|s_i^2-s_j^2|^2\cdot \prod_{i=1}^p s_i^{2(n-p)+1}\cdot\prod_{i=1}^p \exp(-ns_i^2)\,ds_i.
\end{equation*}
\begin{proposition}\label{prop:LDP sym block}
    In the regime $p/n\to 0,$ the averaged singular value distribution $\hat{\mu}_{p,+}^T\in\mathcal{M}(\R_+)$ under the law $W_{p,n}$ satisfies a large deviation principle with speed $pn$ and a good rate function
    \begin{equation*}
	J^+(\mu):=\int (x^2-\log x^2) d\mu(x)-1,\qquad
        \mu\in\mathcal{M}(\R_+).
    \end{equation*}
\end{proposition}
\noindent The minimizer of $J^+(\mu)$ is unique and is given by $\delta_1$ so that we deduce the almost surely convergence of $\hat{\mu}_{p,+}^T$ towards $\delta_1$ as $p\to \infty$. Thus, as a corollary, we get the almost surely convergence of $\hat{\mu}_p^T$ towards the Rademacher distribution $\mathrm{b}$.

\begin{proof}
We mainly follow the standard steps as proving the large deviation principle (see e.g. \cite[Section 2.6.1]{anderson2010introduction}). First, we shall prove the exponential tightness of $W_{p,n}$. Then it suffices to show the large deviation principle for the small ball $B(\mu,\varepsilon),$ where we denote $d$ as the Lipschitz bounded metric in the space $\mathcal{M}(\R_+)$. We organize the proof as follows:

\medskip
\noindent\textbf{Step 1: Exponential tightness.}
We can rewrite the density \eqref{eq:singular density} into:
    \begin{align}\label{eq:sym block rewrite density}
        dW_{p,n}(s_1,\ldots,s_M)=\frac{1}{D_{p,n}}&e^{-p^2\iint_{x\neq y}f(x,y)d\hat{\mu}_{p,+}^T(x)d\hat{\mu}_{p,+}^T(y)}\cdot\\
     &e^{-(n-p)p\int g(x)d\hat{\mu}_{p,+}^T(x)}\prod_{i=1}^pe^{-g(s_i)}ds_i\nonumber ,
    \end{align}
where
    \begin{align*}
        &f(x,y)=\frac{1}{2}\left(x^2+y^2\right)-\log\arrowvert x^2-y^2 \arrowvert,\\
        &g(x)=x^2-\log x^2.
    \end{align*}
    We have two rate functions with two different scales $p(n-p)$ and $p^2$ respectively:
    \begin{align*}
        &J_1^+(\mu):=\iint f(x,y)d\mu(x)d\mu(y)-\inf_{\mu\in\mathcal{M}(\R_+)}\left\{\iint f(x,y)d\mu(x)d\mu(y)\right\},\\
        &J_2^+(\mu):=\int g(x)d\mu(x)-\inf_{\mu\in\mathcal{M}(\R_+)}\left\{\int g(x)d\mu(x)\right\}.
    \end{align*}
By Jensen's inequality, for some constant $C,$
\begin{align*}
    \log D_{p,n}&\geq p\log \int e^{-g(x)}dx
    -p(n-p)\int\left(\int g(x)d\hat{\mu}_{p,+}^T(x)\right)\prod_{i=1}^p\frac{e^{-g(s_i)}ds_i}{\int e^{-g(x)}dx}\\
    &\quad -p^2\int\left(\int_{x\neq y}f(x,y)d\hat{\mu}_{p,+}^T(x)d\hat{\mu}_{p,+}^T(y)\right)\prod_{i=1}^p\frac{e^{-g(s_i)}ds_i}{\int e^{-g(x)}dx}\\
    &\geq -Cpn+p\log \int e^{-g(x)}dx\simeq -Cpn.
\end{align*}
Moreover, notice that there exist some constant $a>0, b>0, c\in\R_+$ such that
\begin{align*}
    &|f(x,y)|\geq a\frac{x^2}{2}+ a\frac{y^2}{2}+ c,\\
    &|g(x)|\geq b\frac{x^2}{2}+c,
\end{align*}
from which one can conclude that for all $K\geq 0,$
\begin{align}\label{ineq:LDP sym block exp tight}
    W_{p,n}\left( \int \frac{x^2}{2}dL_N(x)\geq K\right)\leq e^{-2aKp^2-bKp(n-p)+(C-c)pn}\left( \int e^{-g(x)}d x\right)^p.
\end{align}
Since $\frac{x^2}{2}$ goes to infinity at infinity, the set $\{\mu\in\mathcal{M}(\R_+): \int x^2/2d\mu(x)\geq K \}$ is compact for all $K<+\infty,$ so that we have proved that the law of $\hat{\mu}_{p,+}^T$ under $W_{p,n}$ is exponentially tight.

\medskip
\noindent\textbf{Step 2: Upper bound.}
We set $\bar{W}_{p,n}=D_{p,n}W_{p,n},$ the goal is to prove that for any $\mu\in\mathcal{M}(\R_+)$
\begin{equation*}
    \lim_{\varepsilon\to0}\limsup_{p\to\infty}\frac{1}{pn}\log \bar{W}_{p,n}(d_{BL}(\hat{\mu}_{p,+}^T,\mu)\leq\varepsilon)\leq -\int g(x)d\mu(x).
\end{equation*}
For any $R\geq 0,$ set $f_R(x,y)=f(x,y)\wedge R$ and $ g_R(x)=g(x)\wedge R$. Obviously,
\begin{align*}
    \bar{W}_{p,n}(d_{BL}(\hat{\mu}_{p,+}^T,\mu)\leq \varepsilon)\leq \int_{d_{BL}(\hat{\mu}_{p,+}^T,\mu)\leq\varepsilon} &e^{-p^2\iint_{x\neq y}f_R(x,y)d\hat{\mu}_{p,+}^T(x)d\hat{\mu}_{p,+}^T(y)}\cdot\nonumber\\
    &e^{-(n-p)p\int g_R(x)d\hat{\mu}_{p,+}^T(x)}\prod_{i=1}^pe^{-g(s_i)}ds_i.
\end{align*}
Since $f(x,y)$ is bounded from below and on the set $\{x=y\}, f_R(x,y)=R,$ 
\begin{align*}
    \iint_{x\neq y}f_R(x,y)d\hat{\mu}_{p,+}^T(x)d\hat{\mu}_{p,+}^T(y)&=\iint f_R(x,y)d\hat{\mu}_{p,+}^T(x)d\hat{\mu}_{p,+}^T(y)- R/p\nonumber\\
    &\geq C_f-R/p,
\end{align*}
where $C_f$ is the lower bound of $f$. Combine these two, 
\begin{align}\label{ineq:LDP sym block upperbound}
    &\quad \bar{W}_{p,n}(d_{BL}(\hat{\mu}_{p,+}^T,\mu)\leq \varepsilon)\\
    &\leq e^{-p^2C_f+Rp}\int_{d_{BL}(\hat{\mu}_{p,+}^T,\mu)\leq\varepsilon}
    e^{-(n-p)p\int g_R(x)d\hat{\mu}_{p,+}^T(x)}\prod_{i=1}^pe^{-g(s_i)}ds_i\nonumber.
\end{align}
Note that $g_R(x)$ is bounded continuous, so on $d_{BL}(\hat{\mu}_{p,+}^T,\mu)\leq\varepsilon,$ we have $\int g_R(x)d\mu_{p,+}^P-\int g_R(x)d\mu(x) = o_R(\varepsilon)$ and since $e^{-g}$ is an integrable function, we have for all $R>0,$ 
\begin{equation*}
    \lim_{\varepsilon\to 0}\limsup_{p\to\infty}\frac{1}{pn}\log \bar{W}_{p,n}(d_{BL}(\hat{\mu}_{p,+}^T,\mu)\leq\varepsilon)\leq -\int g_R(x)d\mu(x).
\end{equation*}
Apply monotone convergence theorem and let $R\to \infty,$ we conclude
\begin{equation}\label{ineq:LDP sym block limsup}
    \lim_{\varepsilon\to0}\limsup_{p\to\infty}\frac{1}{pn}\log \bar{W}_{p,n}(d_{BL}(\hat{\mu}_{p,+}^T,\mu)\leq\varepsilon)\leq -\int g(x)d\mu(x).
\end{equation}
The same argument yields that
\begin{equation}\label{ineq:LDP sym block limsupconstant}
    \limsup_{p\to \infty}\frac{1}{pn}\log D_{p,n}\leq -\inf_{\mu\in\mathcal{M}(\R_+)}\left\{\int g(x)d\mu(x)\right\}.
\end{equation}

\medskip
\noindent\textbf{Step 3: Lower bound.}
This part aims to prove the following
\begin{equation}\label{ineq:LDP sym block liminf}
    \lim_{\varepsilon\to0}\liminf_{p\to\infty}\frac{1}{pn}\log \bar{W}_{p,n}(d_{BL}(\hat{\mu}_{p,+}^T,\mu)\leq\varepsilon)\geq -\int g(x)d\mu(x).
\end{equation}
Without loss of generality, we can assume that $\int g(x)d\mu(x)<\infty,$ which implies that the distribution function of $\mu$ is continuous near $0$. So if we set $K_{\delta}^L=[\delta,L]$ with $0<\delta<L,$ it suffices to consider the probability measure supported in $K_{\delta}^L$. Indeed, for any $\mu\in\mathcal{M}(\R_+),$ the truncated probability measure $\mu_{\delta}^L=(\mu(K_{\delta}^L))^{-1}\1_{ K_{\delta}^L}d\mu(x)$ converges weakly to $\mu$ as $\delta\to 0, L\to \infty$. Since $g(x)$ is bounded from below, again by the monotone convergence theorem,
\begin{equation*}
    \lim_{\substack{\delta\to 0\\ L\to \infty}}\int g(x)d\mu_{\delta}^L(x)=\int g(x)d\mu(x).
\end{equation*}
Moreover, it is enough to prove the case where $\mu$ has no atom since any distribution function can be approximated by a continuous one. For any $\varepsilon>0,$ there exists $\mu_{\varepsilon}\in\mathcal{M}(\R_+)$ with no atoms and support contained in $K_{\delta}^L$ such that $d_{BL}(\mu_{\varepsilon},\mu)<\varepsilon/2$. By the triangle inequality, for any $\eta<\varepsilon/2,$ we have
\begin{equation*}
    \bar{W}_{p,n}(d_{BL}(\hat{\mu}_{p,+}^T,\mu)\leq \varepsilon)\geq \bar{W}_{p,n}(d_{BL}(\hat{\mu}_{p,+}^T,\mu_{\varepsilon})\leq \eta).
\end{equation*}
Let $p\to \infty, \eta\to 0, \varepsilon\to 0$ and since $g(x)$ is bounded continuous in $K_{\delta}^L,$ 
\begin{align*}
    \lim_{\varepsilon\to 0}\liminf_{p\to \infty}\frac{1}{pn}&\log \bar{W}_{p,n}(d_{BL}(\hat{\mu}_{p,+}^T,\mu)\leq \varepsilon)\\
    &\geq \lim_{\varepsilon\to 0}\lim_{\eta\to 0}\liminf_{p\to \infty}\frac{1}{pn}\log \bar{W}_{p,n}(d_{BL}(\hat{\mu}_{p,+}^T,\mu_{\varepsilon})\leq \eta)\\
    &\geq \lim_{\varepsilon\to 0} -\int g(x)d\mu_{\varepsilon}(x)= -\int g(x)d\mu(x).
\end{align*}
In this way, we can obtain the lower bound for general $\mu\in\mathcal{M}(\R_+)$. Therefore, it remains to show that if $\mu\in\mathcal{M}(\R_+)$ has no atoms and is supported in $K_{\delta}^L,$ 
\begin{equation*}
    \lim_{\varepsilon\to0}\liminf_{p\to\infty}\frac{1}{pn}\log \bar{W}_{p,n}(d_{BL}(\hat{\mu}_{p,+}^T,\mu)\leq\varepsilon)\geq -\int g(x)d\mu(x).
\end{equation*}
The idea is to localize the singular values $(s_i)_{1\leq i\leq p}$ in small sets and to take advantage of the speed of $p(n-p)$ to neglect the small volume of these sets. To do so, first remark that for any $\mu\in\mathcal{M}(\R_+)$ with no atoms, if we set
\begin{align*}
    &x^{1,p}=\inf\left\{x: \mu((-\infty,x])\geq \frac{1}{p+1}\right\},\\
    &x^{i+1,p}=\inf\left\{x\geq x^{i,p}: \mu((x^{i,p},x])\geq\frac{1}{p+1}\right\},\quad 1\leq i\leq p-1,
\end{align*}
for any real number $\eta,$ there exists a positive integer $p(\eta)$ such that, for any $p$ larger than $p(\eta),$
\begin{equation*}
    d_{BL}\left( \mu, \frac{1}{p}\sum_{i=1}^p\delta_{x_i}\right)<\eta.
\end{equation*}
In particular, for $p\geq p(\frac{\varepsilon}{2}),$
\begin{align*}
    \left\{(s_i)\in\R_+^p:  |s_i-x^{i,p}|<\frac{\varepsilon}{2}, 1\leq i\leq p \right\}\subset\{(s_i)\in\R_+^p: d_{BL}(\hat{\mu}_{p,+}^T,\nu)<\varepsilon\}. 
\end{align*}
Now we take the associated division $\{x^{i,p},i=1,\ldots,p\}$ for $\mu\in\mathcal{M}(\R_+)$ so that we have the lower bound:
\begin{align}\label{ineq:LDP sym block lowerbound}
    &\bar{W}_{p,n}(d_{BL}(\hat{\mu}_{p,+}^T,\mu)\leq \varepsilon)\\
    &\geq \int_{\cap_i\{|s_i-x^{i,p}|<\frac{\varepsilon}{2}\}}\prod_{1\leq i<j\leq p}|s_i^2-s_j^2|^2\cdot \prod_{i=1}^p s_i^{2(n-p)+1}\cdot\prod_{i=1}^p \exp(-ns_i^2)\1_{\cap\{s_i\geq 0\}}\prod_{i=1}^pds_i\nonumber\\
    &= \int_{\cap_i\{|s_i|<\frac{\varepsilon}{2}\}}\prod_{i=1}^p|s_i+x^{i,p}|^{2(n-p)+1}\nonumber\\
    &\quad \cdot \prod_{i<j}|(s_i+x^{i,p})^2-(s_j+x^{j,p})^2|^2  \cdot \prod_{i=1}^pe^{-n(s_i+x^{i,p})^2}ds_i\nonumber\\
    &\geq \left(\prod_{i=1}^p \left(x^{i,p}\right)^{2(n-p)+1}\cdot e^{-n\left(x^{i,p}\right)^2}\right) \cdot \int_{\substack{{\cap_i\{|s_i|<\frac{\varepsilon}{2}\}}\\  s_i<s_{i+1}}} \prod_{i<j}|(s_i+x^{i,p})^2-(s_j+x^{j,p})^2|^2\cdot\nonumber\\
    &\quad\quad e^{-(p-1/2)\sum_{i=1}^p[(x^{i,p}+s_i)^2-(x^{i,p})^2]}\cdot e^{-(n-p+1/2)\sum_{i=1}^p[g(x^{i,p}+s_i)-g(x^{i,p})]}\prod_{i=1}^p ds_i\nonumber\\
    &\geq\left(\prod_{i=1}^p \left(x^{i,p}\right)^{2(n-p)+1}\cdot e^{-\frac{n}{2}\left(x^{i,p}\right)^2}\right)\cdot \left(\frac{\delta}{2}\right)^{p(p-1)/2}\times \nonumber\\
    &\quad \Bigg(\int_{\substack{{\cap_i\{|s_i|<\frac{\varepsilon}{2}\}}\\s_i<s_{i+1}}} \prod_{i<j} |s_i-s_j|^2 \cdot e^{-(p-1/2)\sum_{i=1}^p[(x^{i,p}+s_i)^2-(x^{i,p})^2]}\nonumber\\
    &\quad \cdot e^{-(n-p+1/2)\sum_{i=1}^p[g(x^{i,p}+s_i)-g(x^{i,p})]}\prod_{i=1}^p ds_i\Bigg)\nonumber\\
    &:= \bar{W}_{p,n}^1\times \bar{W}_{p,n}^2,\nonumber
\end{align}
here the first equality is by the change of variables $s_i-x^{i,p}= s'_i.$ The third inequality is due to the fact that
\begin{align*}
    &|x^{i,p}-x^{j,p}+s_i-s_j|\geq |x^{i,p}-x^{j,p}|\vee |s_i-s_j|\geq |s_i-s_j |,\\
    &|x^{i,p}+x^{j,p}+s_i+s_j|\geq 2|s_1+x_1|\geq \frac{\delta}{2},
\end{align*}
provided $0<s_i\leq s_j,\ \delta<x^{i,p}\leq x^{j,p},$ for all $i<j$ and $\varepsilon>0$ small enough.

To deal with the term $\bar{W}_{p,n}^1,$ by the choice of $x^{i,p},$
\begin{align*}
    &\lim_{p\to \infty}\frac{1}{p}\sum_{i=1}^p(x^{i,p})^2=\int x^2d\mu(x),\nonumber\\
    &\lim_{p\to \infty}\frac{1}{p}\sum_{i=1}^p\log(x^{i,p})^2=\int \log x^2d\mu(x).
\end{align*}
Hence
\begin{align}\label{eq:LDP sym block limQ1bar}
    \lim_{p\to \infty}\frac{1}{pn}\log \bar{W}_{p,n}^1 &=-\lim_{p\to \infty}\frac{n-p+1/2}{pn}\sum_{i=1}^p\log (x^{i,p})^2-\lim_{p\to \infty}\frac{1}{p}\sum_{i=1}^p(x^{i,p})^2\nonumber\\
    &\quad -\lim_{p\to\infty}\frac{p(p-1)}{2pn}\log \frac{\delta}{2}\\
    &=-\int g(x)d\mu(x)\nonumber.
\end{align}

To estimate $\bar{W}_{p,n}^2,$ note that since we assumed that $\mu$ have compact support $K_L^{\delta},$ $(x^{i,p},1\leq i\leq p)$ are uniformly bounded and so by the uniform continuity of $x^2$ and $\log x^2$ in $K_L^{\delta},$
\begin{align}\label{eq:LDP sym block uniformcontinuous}
    &\lim_{\varepsilon\to 0}\sup_{p\in\mathbb{N}}\sup_{1\leq i\leq p}\sup_{|x|\leq \varepsilon}|(x^{i,p}+x)^2-(x^{i,p})^2|=0,\\
    &\lim_{\varepsilon\to 0}\sup_{p\in\mathbb{N}}\sup_{1\leq i\leq p}\sup_{|x|\leq \varepsilon}|\log(x^{i,p}+x)^2-\log (x^{i,p})^2|=0.\nonumber
\end{align}
Moreover, by symmetry and Jensen's inequality,
\begin{align}
    &\quad \log\int_{\substack{{\cap_i\{|s_i|<\frac{\varepsilon}{2}\}}\\s_1<s_2<\cdots<s_M}} \prod_{1\leq i<j\leq p} |s_i-s_{j}|^2 \prod_{i=1}^p ds_i\nonumber\\
    &= \log \frac{1}{p!}\int_{\cap_i\{|s_i|<\frac{\varepsilon}{2}\}}\prod_{1\leq i<j\leq p}|s_i-s_j|^2\prod_{i=1}^pds_i\nonumber\\
    &\geq \log\frac{\varepsilon^p}{p!}+\frac{1}{\varepsilon^p}\int_{\cap_i\{|s_i|<\frac{\varepsilon}{2}\}}\sum_{1\leq i<j\leq p} \log|s_i-s_{j}|^2 \prod_{i=1}^p ds_i\nonumber\\
    &\geq \log\frac{\varepsilon^p}{p!}+\sum_{1\leq i<j\leq p}\frac{1}{\varepsilon^2}\int_{[-\frac{\varepsilon}{2},\frac{\varepsilon}{2}]^2} \log|s_i-s_j|^2 ds_ids_j\nonumber
\end{align}
By the change of variables: $s_j= y,\ s_i-s_j= x,$ we get
\begin{equation*}
    \int_{[-\frac{\varepsilon}{2},\frac{\varepsilon}{2}]^2} \log|s_i-s_j|^2 ds_ids_j= \int_{-\varepsilon/2}^{\varepsilon/2}dy\int_{-\varepsilon}^{\varepsilon} \log x^2dx=4\varepsilon^2(\log \varepsilon-1)
\end{equation*}
Note that by Stirling's approximation, $p!\sim \sqrt{2\pi p}\left(p/e\right)^p,$ so
\begin{equation*}
    \liminf_{p\to \infty}\frac{1}{pn}\log\int_{\substack{{\cap_i\{|s_i|<\frac{\varepsilon}{2}\}}\\s_1<s_2<\cdots<s_M}} \prod_{1\leq i<j\leq p} |s_i-s_{j}|^2 \prod_{i=1}^p ds_i\geq 0
\end{equation*}
and together with \eqref{eq:LDP sym block uniformcontinuous}, we have
\begin{equation}\label{ineq:LDP sym block liminfQ2bar}
    \lim_{\varepsilon\to 0}\liminf_{p\to \infty}\frac{1}{pn}\log \bar{W}_{p,n}^2\geq 0.
\end{equation}
Combine \eqref{eq:LDP sym block limQ1bar} and \eqref{ineq:LDP sym block liminfQ2bar}, we get the desired lower bound \eqref{ineq:LDP sym block liminf}.

\medskip
\noindent\textbf{Step 4: Combining Steps 1-3 to complete the proof.}
For any $\mu\in\mathcal{M}(\R_+),$ by \eqref{ineq:LDP sym block liminf}
\begin{align*}
    \liminf_{p\to \infty}\frac{1}{pn}\log D_{p,n}&\geq \lim_{\varepsilon\to 0}\liminf_{p\to \infty}\frac{1}{pn}\log \bar{W}_{p,n}(d_{BL}(\hat{\mu}_{p,+}^T,\mu)\leq \varepsilon)\\
    &\geq -\int g(x)d\mu(x),
\end{align*}
take the supremum over $\mu$ on the right-hand side of the inequality and combine it with \eqref{ineq:LDP sym block limsupconstant},
\begin{equation}\label{eq:LDP sym block limconstant}
    \lim_{p\to \infty}\frac{1}{pn}\log D_{p,n}=-\inf_{\mu\in\mathcal{M}(\R_+)}\left\{\int g(x)d\mu(x)\right\}.
\end{equation}
Thus, \eqref{ineq:LDP sym block limsup}, \eqref{ineq:LDP sym block liminf} and \eqref{eq:LDP sym block limconstant} imply the weak large deviation principle i.e.
\begin{align*}
    &\quad\lim_{\varepsilon\to 0}\liminf_{p\to \infty}\frac{1}{pn}\log W_{p,n}(d_{BL}(\hat{\mu}_{p,+}^T,\mu)\leq \varepsilon)\nonumber\\
    &=\lim_{\varepsilon\to 0}\limsup_{p\to \infty}\frac{1}{pn}\log W_{p,n}(d_{BL}(\hat{\mu}_{p,+}^T,\mu)\leq \varepsilon)\nonumber\\
    &=-\left(\int g(x)d\mu(x)-\inf_{\mu\in\mathcal{M}(\R_+)}\left\{\int g(x)d\mu(x)\right\}\right)=-J_2^+(\mu).
\end{align*}
This, together with the exponential tightness \eqref{ineq:LDP sym block exp tight}, completes the proof. 
\end{proof}

\subsection{Link with Wishart matrix}
Applying a similar argument, we obtain a large deviation principle for the
empirical spectral measure $\tilde{\mu}_p=p^{-1}\sum_{i=1}^p\delta_{\lambda_i(W(p,n))}$ of the Wishart matrix $W(p,n)$ in the regime $p/n\to0$. Recall that the
eigenvalue density of $W(p,n)$ is proportional to
\[
\prod_{1\le i<j\le p}|\lambda_i-\lambda_j|^2
\prod_{i=1}^p\lambda_i^{\,n-p}
\prod_{i=1}^p e^{-n\lambda_i}\,d\lambda_i .
\]

\begin{corollary}
In the regime $p/n\to0$, the empirical spectral measure $\tilde{\mu}_p$
satisfies a large deviation principle with speed $pn$ and good rate
function
\[
J^+(\mu)
=
\int (x-\log x)\,d\mu(x)-1,
\qquad
\mu\in\mathcal M(\mathbb R_+).
\]
The unique minimizer of $J^+$ is given by $\delta_1$.
\end{corollary}

This result can be viewed as a complement to the results of Hiai and Petz
\cite{hiai1998eigenvalue}, where a large deviation principle for the same
model was established in the regime $p/n\to\gamma\in(0,1].$
More precisely, using the speed $pn$, the corresponding rate function on
$\mathcal M(\mathbb R_+)$ is given by
\begin{equation}\label{eq:rate function wishart}
   J_\gamma^+(\mu)
=
\gamma\left(
\int x\,d\mu(x)-\Sigma(\mu)
\right)
+
(1-\gamma)
\int (x-\log x)\,d\mu(x)
-c_\gamma ,    
\end{equation}
where
\[
c_\gamma
=
\inf_{\mu\in\mathcal M(\mathbb R_+)}
\left\{
\gamma\left(
\int x\,d\mu(x)-\Sigma(\mu)
\right)
+
(1-\gamma)
\int (x-\log x)\,d\mu(x)
\right\}.
\]

The unique minimizer of $J_\gamma^+$ is the Marchenko--Pastur distribution
$\nu_\gamma$ with density
\begin{equation}\label{eq:MP density}
    d\nu_\gamma(t)
=
\frac{\sqrt{
(t-(1-\sqrt{\gamma})^2)
((1+\sqrt{\gamma})^2-t)
}}{2\pi\gamma t}\,\1_{[(1-\sqrt{\gamma})^2,(1+\sqrt{\gamma})^2]}(t)\,dt .
\end{equation}
Thus, the empirical spectral measure $\tilde{\mu}_p$ converges almost surely to
$\nu_\gamma$. Moreover, we observe that as $\gamma\to0$, $\nu_\gamma\to\delta_1.$ in distribution. To understand the fluctuations around this limiting behavior, Jiang
\cite{jiang2015approximation} first studied the shifted and rescaled
eigenvalues
\[
\sigma_i=\frac{\lambda_i-\beta n}{\sqrt{2\beta pn}},
\]
where $(\lambda_1,\ldots,\lambda_p)$ follows the $\beta$-Laguerre ensemble.
He proved a weak large deviation principle for the rescaled empirical spectral measure $p^{-1}\sum_{i=1}^p\delta_{\sigma_i}$
under the assumption $p^2/n\to0.$ More recently, Ma \cite{ma2023unified} considered the same rescaled empirical spectral measure
with a different normalization
\[
\sigma_i'=\frac{\lambda_i-\beta n}{2\beta\sqrt{pn}},
\]
where the parameter $\beta$ is allowed to vary with $n$. Under the assumptions $p/n\to\gamma\in[0,1],\ 
(\log p)/\beta p\to0,$ he established a full large deviation principle with speed $\beta p^2$ for $\frac1p\sum_{i=1}^p\delta_{\sigma_i'}.$
with speed $\beta p^2$ with a good rate function. In particular, in the regime $p/n\to0$, the rate function coincides with the one of GUE \eqref{eq:rate function GUE} and thus the rescaled empirical spectral measure converges almost surely to the semicircle law.

\section{LDP for the Conditioned GUE model}
\subsection{Proving Theorem \ref{thm:LDP con GUE}}
We sketch the proof of Theorem \ref{thm:LDP con GUE}. Recall that the density function \eqref{eq:con GUE density} of the eigenvalues $(\lambda_1,\ldots,\lambda_M)$ is given by
\begin{equation*}
    d Q_{M,N}(\lambda_1,\ldots,\lambda_M)=\frac{1}{Z_{M,N}}\prod_{i=1}^M\lambda_i^{2(N-M)}\prod_{1\leq i<j\leq M}|\lambda_i-\lambda_j|^2 \prod_{i=1}^M e^{-N\frac{\lambda_i^2}{2}}d\lambda_i.
\end{equation*}
As before, we can write this density in the following form:
\begin{align*}
    dQ_{M,N}(\lambda_1,\ldots,\lambda_M)=\frac{1}{Z_{M,N}}&e^{-M^2\iint_{x\neq y}F(x,y)dL_M(x)dL_M(y)}\cdot\\
    &e^{-(N-M)M\int G(x)dL_M(x)}\prod_{i=1}^Me^{-\frac{\lambda_i^2}{2}}d\lambda_i,
\end{align*}
where
\begin{align*}
    &F(x,y)=\frac{1}{2}\left(\frac{1}{2}x^2+\frac{1}{2}y^2\right)-\log\arrowvert x-y \arrowvert,\\
    &G(x)=\frac{1}{2}x^2-\log x^2.
\end{align*}
We have two rate functions with two different scales $M(N-M)$ and $M^2$ respectively:
\begin{align*}
    &I_1(\mu):=\iint F(x,y)d\mu(x)d\mu(y)-\inf_{\mu\in\mathcal{M}(\mathbb{R})}\left\{\iint F(x,y)d\mu(x)d\mu(y)\right\},\\
    &I_2(\mu):=\int G(x)d\mu(x)-\inf_{\mu\in\mathcal{M}(\mathbb{R})}\left\{\int G(x)d\mu(x)\right\}.
\end{align*}
with different scales $M(N-M)$ and $M^2$ respectively. Similarly, one can verify that the arguments in the proof of Theorem 
\ref{thm:LDP sym block} remain valid for this model. Therefore, the speed 
of the large deviation principle is expected to be $MN$, with the 
corresponding rate function denoted by $I_2$. However, we note that the 
minimizers of the rate function $I_2$ are not unique. They are given by the 
set of atomic probability measures supported on $\{\pm\sqrt{2}\}$, namely,
\[
\mathcal{M}_0
=
\left\{
\mu_p:=p\delta_{\sqrt{2}}+(1-p)\delta_{-\sqrt{2}}
:\ p\in[0,1]
\right\}.
\]
Combining this observation with Theorem \ref{thm:LDP con GUE}, we obtain
\begin{corollary}\label{cor:con GUE a.s.}
Almost surely, $d_{BL}(L_M,\mathcal{M}_0)
\stackrel{M\to\infty}{\longrightarrow}0.$
Moreover,
\[
\mathbb{E}[L_M]\xrightarrow{\mathcal{L}}
\frac12\delta_{-\sqrt2}+\frac12\delta_{\sqrt2}.
\]
\end{corollary}

\begin{proof}
Fix $\varepsilon>0$ and denote by $\mathcal{M}_0^\varepsilon$ the 
$\varepsilon$-neighborhood of $\mathcal{M}_0$. Since $\mathcal{M}_0$ is the 
set of minimizers of the rate function, there exists
$\gamma_\varepsilon>0$ such that
\[
\inf_{\mu\notin\mathcal{M}_0^\varepsilon}I(\mu)
\geq \gamma_\varepsilon .
\]
By the large deviation principle, for sufficiently large $N$,
\[
\frac1{NM}
\log\mathbb P(L_M\notin\mathcal M_0^\varepsilon)
\leq
-\frac12
\inf_{\mu\notin\mathcal M_0^\varepsilon}I(\mu)
\leq
-\frac{\gamma_\varepsilon}{2}.
\]
Therefore,
\[
\sum_{M=1}^{\infty}
\mathbb P(L_M\notin\mathcal M_0^\varepsilon)<\infty .
\]
By the Borel--Cantelli lemma, $\mathbb P(L_M\notin\mathcal M_0^\varepsilon\ {\rm i.o.})=0.$ Since this holds for every $\varepsilon>0$, the first claim follows.

For the second claim, for any odd integer $k$, by symmetry,
\[
\mathbb E\left[\int x^k\,dL_M(x)\right]
=0
=
\int x^k\,d\mu_{1/2}(x).
\]
For even $k$, observe that for any $p,q\in[0,1]$,
\[
\int x^k\,d\mu_p(x)
=
\int x^k\,d\mu_q(x).
\]
Hence, using the convergence of $L_M$ towards the minimizer set
$\mathcal M_0$, we obtain that for every $k\ge1$,
\[
\mathbb E\left[\int x^k\,dL_M(x)\right]
\longrightarrow
\int x^k\,d\mu_{1/2}(x),
\]
which completes the proof.
\end{proof}

\subsection{Another LDP for the scaled empirical spectral measure}
We aim to show that $L_M$ converges weakly almost surely to $\mu_{\frac12}$. 
Motivated by the works of Jiang \cite{jiang2015approximation} and Ma 
\cite{ma2023unified} discussed in the previous section, our strategy is to 
zoom the eigenvalues around the two points $\pm\sqrt{2}$ with an 
appropriate scaling factor and derive another large deviation principle.

We first introduce some notations. For a given configuration $(\lambda_i)_{1\leq i\leq M}$ distributed according to the law in \eqref{eq:con GUE density}, define
\[
M_0=\#\{i:\lambda_i\geq 0\}.
\]
By symmetry, we consider the set $\Delta=\{(\lambda_1,\ldots,\lambda_M)\in \mathbb{R}^M,\,\lambda_1\ge \lambda_2\ge \cdots\ge\lambda_M\}$, and define the rescaled variables
\begin{align*}
	\alpha_i&=\Theta_M^{-1}(\lambda_i-\sqrt{2}),\qquad i=1,\ldots,M_0,\ (\lambda_i\ge 0),\\
	\beta_i&=\Theta_M^{-1}(\lambda_j+\sqrt{2}),\qquad i=M_0+1,\ldots,M,\ (\lambda_i<0),
\end{align*}
where $\Theta_M$ is a scaling parameter satisfying $\Theta_M\to0$ as $M\to\infty$.

Recall that the density of $(\lambda_1,\ldots,\lambda_M)\in \Delta$ is given by
\begin{equation*}
d Q_{M,N}(\lambda_1,\ldots,\lambda_M)
=
\frac{M!}{Z_{M,N}}
\prod_{i=1}^M\lambda_i^{2(N-M)}
\prod_{1\leq i<j\leq M}|\lambda_i-\lambda_j|^2
\prod_{i=1}^M e^{-N\frac{\lambda_i^2}{2}}d\lambda_i.
\end{equation*}
By the elementary decomposition, for any Borel set $A\in \mathbb{R}^M$,
\begin{equation*}
   \mathbb{P}((\lambda_1,\ldots,\lambda_M)\in A) = \sum_{k=0}^{M} \mathbb{P}((\alpha_1,\ldots,\alpha_k,\beta_{k+1},\ldots,\beta_M)\in A_k),
\end{equation*}
where $A_k$ is the image of $A\cap\{M_0=k\}$ with respect to the maps $\lambda_i\mapsto \alpha_i,\,i=1,\ldots,k;\,\lambda_i\mapsto\beta_i,\,i=k+1,\ldots,M$. Thus
\begin{align}\label{eq:con GUE density alpha beta}
    &\quad \mathbb{P}((\lambda_1,\ldots,\lambda_M)\in A)\nonumber\\
    &=\frac{1}{\tilde{Z}_{M,N}}
    \sum_{k=0}^M\int_{A_k}
    \prod_{i=1}^k
    \left[
    (\sqrt{2}+\Theta_M\alpha_i)^{2(N-M)}
    \exp\left(-\frac{N}{2}(\sqrt{2}+\Theta_M\alpha_i)^2\right)
    \right]
    \nonumber\\
    &\quad \times
    \prod_{j=k+1}^{M}
    \left[
    (-\sqrt{2}+\Theta_M\beta_j)^{2(N-M)}
    \exp\left(-\frac{N}{2}(-\sqrt{2}+\Theta_M\beta_j)^2\right)
    \right]
    \nonumber\\
    &\quad \times
    \Theta_M^{-2k(M-k)}
    \prod_{1\le i<i'\le k}|\alpha_i-\alpha_{i'}|^2
    \prod_{k+1\le j<j'\le M}|\beta_j-\beta_{j'}|^2
    \nonumber\\
    &\quad \times
    \prod_{\substack{1\leq i\leq k\\k+1\leq j\leq M}}
    |2\sqrt{2}+\Theta_M(\alpha_i-\beta_j)|^2
    \cdot
     d\alpha_id\beta_j ,
\end{align}
where $\tilde{Z}_{M,N}=Z_{M,N}\Theta_M^{-M^2-M} (M!)^{-1}$. Define the associated empirical spectral measure of $\{\alpha_i\}_{i=1}^{M_0}$ and $\{\beta_j\}_{j=M_0+1}^M$ as follows:
\[
\mu_{\alpha,M}
=
\frac1M\sum_{i=1}^{M_0}\delta_{\alpha_i},
\qquad
\mu_{\beta,M}
=
\frac1M\sum_{j=M_0+1}^M\delta_{-\beta_j}.
\]
We denote by $\mathcal{M}_{\leq1}(\mathbb{R})$ the space of finite positive measures on $\mathbb{R}$ with total mass at most $1$. We equip $\mathcal{M}_{\leq1}(\mathbb{R})$ with the weak topology induced by the bounded-Lipschitz metric $d_{BL}$. For a finite positive measure $\mu$, we write
\[
\operatorname{mass}(\mu):=\int_{\mathbb{R}}1\,d\mu
\]
for its total mass. Now consider the subspace
\begin{equation*}
    \mathcal{M}
    =
    \{(\mu,\nu)\in \mathcal{M}_{\leq1}(\mathbb{R})\times
    \mathcal{M}_{\leq1}(\mathbb{R}):
    \operatorname{mass}(\mu)+\operatorname{mass}(\nu)=1\},
\end{equation*}
which is endowed with the topology inherited from 
$\mathcal{M}_{\leq1}(\mathbb{R})\times\mathcal{M}_{\leq1}(\mathbb{R})$
equipped with the metric $d_{BL}\oplus d_{BL}$. Our goal is to establish a large deviation principle for the pair
$(\mu_{\alpha,M},\mu_{\beta,M})\in\mathcal{M}$.

\medskip
\noindent\textbf{Heuristics.} To identify the correct large deviation speed, we first rewrite the density \eqref{eq:con GUE density alpha beta} in terms of $\mu_{\alpha,M},\mu_{\beta,M}$, for the first term, we divide into two parts
\begin{align*}
	&\quad C_{\alpha,N}\prod_{i=1}^k \left[(\sqrt{2}+\Theta_M\alpha_i)^{2(N-M)}\exp\left(-\frac{N-M}{2}(\sqrt{2}+\Theta_M\alpha_i)^2\right)\right]\\
	&=\exp\left[-(N-M)M\int \left(\frac{\Theta_M}{\sqrt{2}} x\right)^2+2\frac{\Theta_M}{\sqrt{2}} x-2\log \left(1+\frac{\Theta_M}{\sqrt{2}} x\right)d\mu_{\alpha,M}(x)\right]\\
	&:= \exp\left[-(N-M)MI_1^{\alpha}\right],
\end{align*}
where $C_{\alpha,N}=\exp\left[-(N-M)M(1-\log2)\text{mass}(\mu_{\alpha,M})\right].$
The same computation applies to $\mu_{\beta,M}$ and we denote by $I_1^{\beta}$. Note that if we combine $I_1^\alpha$ with $I_1^\beta$, the normalized factor $C_{\alpha,N}\cdot C_{\beta,N}= \exp\left[(N-M)M(1-\log2)\right]$ only depends on $M$, and we may just put it into the normalized constant $\tilde{Z}_{M,N}.$  The second part is given by 
\begin{align*}
	\prod_{i=1}^k\exp\left(-\frac{M}{2}(\sqrt{2}+\Theta_M\alpha_i)^2\right)&=\exp\left[-\frac{M^2}{2}\int (\sqrt{2}+\Theta_Mx)^2d\mu_{\alpha,M}(x)\right]\\
	&:=\exp[M^2I_2^{\alpha}].
\end{align*} 
and the same definition works for $I_2^\beta$.

For the self-interacting term,
\begin{align*}
	\prod_{i<i'} |\alpha_i-\alpha_{i'}|^2 &=\exp\left[ M^2\iint_{x\neq y} \log |x-y|d\mu_{\alpha,M}(x)d\mu_{\alpha,M}(y)\right]\\
    &:=\exp\left[M^2I_2^{\alpha,\alpha}\right],
\end{align*}
and we define $I_2^{\beta,\beta}$ in the same way. Also,
\begin{align*}
	\prod_{\substack{1\leq i\leq k\\k+1\leq j\leq M}}
    &|2\sqrt{2}+\Theta_M(\alpha_i-\beta_j)|^2\\
    &=\exp\left[M^2 \iint \log(2\sqrt{2}+\Theta_M(x+y))^2d\mu_{\alpha,M}(x)d\mu_{\beta,M}(y)\right]\\
	&:=\exp[M^2I_2^{\alpha,\beta}].
\end{align*}
Now we are able to rearrange the density \eqref{eq:con GUE density alpha beta} after combining constants into the partition function $\tilde{Z}_{M,N}$ (without loss of generality we still denote it as $\tilde{Z}_{M,N}$):
\begin{align*}
	&\frac{1}{\tilde{Z}_{M,N}}\sum_{k=0}^M \exp\left[-M^2(I_2^{\alpha}+I_2^{\beta}+I_2^{\alpha,\alpha}+I_2^{\beta,\beta}+I_2^{\alpha,\beta})\right]\cdot\\
	&\quad \exp\left[k(M-k)\log \Theta_M^{-2}\right]\exp \left[-M(N-M)(I_1^{\alpha}+I_1^{\beta})\right]\cdot
     \prod_{\substack{1\leq i\leq k\\k+1\leq j\leq M}}d\alpha_id\beta_j\nonumber.
\end{align*} 
At the heuristic level of large deviation analysis, the terms
$I_2^{\alpha}$, $I_2^{\beta}$, $I_2^{\alpha,\alpha}$, 
$I_2^{\beta,\beta}$, and $I_2^{\alpha,\beta}$ are of order $M^2$.
Moreover, since $\Theta_M$ is small, a Taylor expansion of $\log(1+x)$ gives
\begin{equation*}
    \left(\frac{\Theta_M}{\sqrt{2}}x\right)^2
    +2\frac{\Theta_M}{\sqrt{2}}x
    -2\log\left(1+\frac{\Theta_M}{\sqrt{2}}x\right)
    \sim
    \Theta_M^2x^2.
\end{equation*}
Hence, the contribution of the corresponding term in the density is roughly of order
\[
M(N-M)(I_1^\alpha+I_1^\beta)
\sim NM\Theta_M^2.
\]

On the other hand, the term $k(M-k)\log\Theta_M^{-2}$ between the two groups contributes
$$
\exp\{-M^2\log \Theta_M^{-2}\iint 1 \,d\mu_{\alpha,M}(x)d\mu_{\beta,M}(y)\}
$$
of order $M^2\log\Theta_M^{-2}$. Therefore, in order for both contributions to appear at the same scale, we choose
$\Theta_M$ according to the balance relation $M^2\log\Theta_M^{-2}=NM\Theta_M^2.$ Since $\Theta_M\to0$, we have $M^2\log\Theta_M^{-2}\gg M^2.$ Thus, similarly to the previous section where a large deviation principle with speed $NM$ suppresses terms of order $M^2$, the terms $I_2^\alpha$, $I_2^\beta$, $I_2^{\alpha,\alpha}$, $I_2^{\beta,\beta}$, and $I_2^{\alpha,\beta}$ are negligible at the present scale.

Consequently, we expect that the pair $(\mu_{\alpha,M},\mu_{\beta,M})\in\mathcal{M}$ satisfies a large deviation principle with speed $\Lambda_M:=M^2\log\Theta_M^{-2},$
and a rate function given by (up to an additive constant arising from the normalization)
\begin{align*}
I_{\alpha,\beta}(\mu_\alpha,\mu_\beta)
\triangleq
\sum_{\sigma=\alpha,\beta}
\int x^2\,d\mu_\sigma(x)
-
\iint 1\,d\mu_\alpha(x)d\mu_\beta(y).
\end{align*}
This heuristic is made precise in the following theorem.

\begin{theorem}\label{thm:LDP rescaled con GUE}
Suppose that $\Theta_M$ satisfies $M^2\log\Theta_M^{-2}=NM\Theta_M^2.$ Then, under the law $\tilde{Q}_{M,N}$, the pair
$(\mu_{\alpha,M},\mu_{\beta,M})\in\mathcal{M}$ satisfies a large deviation principle with speed $\Lambda_M:=M^2\log\Theta_M^{-2},$
and good rate function
$$
\bar{I}_{\alpha,\beta}(\mu_\alpha,\mu_\beta)
:=
I_{\alpha,\beta}(\mu_\alpha,\mu_\beta)
-
\inf_{(\mu,\nu)\in\mathcal{M}}
I_{\alpha,\beta}(\mu,\nu).
$$
\end{theorem}

\begin{proof}
    The proof is divided into three parts, following the standard approach. We first establish the \textbf{exponential tightness}. Since the argument is essentially similar to the one \eqref{ineq:LDP sym block exp tight} in the previous section, we omit the details.

\medskip
\noindent\textbf{Upper bound.} For $(\mu_{\alpha},\mu_{\beta})\in\mathcal{M}$, denote
\[
A_{\alpha,\beta,M}^{\varepsilon}
=\{(\lambda_i)\in \mathbb{R}^M: d_{BL}((\mu_{\alpha,M},\mu_{\beta,M}),(\mu_\alpha,\mu_\beta))\le \varepsilon\}
,
\]
hence, $\mathbb{P}(d_{BL}((\mu_{\alpha,M},\mu_{\beta,M}),(\mu_\alpha,\mu_\beta))\le \varepsilon)=Q_{M,N}(A_{\alpha,\beta,M}^\varepsilon)$.

Denote $d\tilde{Q}_{M,N} = \tilde{Z}_{M,N}dQ_{M,N},$ we aim to prove that for any
$(\mu_{\alpha},\mu_{\beta})\in\mathcal{M}$,
\begin{align}\label{ineq:LDP con GUE alpha beta limsup}
    \lim_{\varepsilon\to0}\limsup_{M\to\infty}
    \frac{1}{\Lambda_M}
    \log\tilde{Q}_{M,N}
    (A_{\alpha,\beta,M}^{\varepsilon})
    \leq
    -I_{\alpha,\beta}(\mu_{\alpha},\mu_{\beta}).
\end{align}

The strategy is analogous to the proof of the large deviation principle in the previous section. We first rewrite the terms $I_i$ into more familiar forms.
\begin{align*}
    I_2^{\alpha,\alpha}+I_2^{\beta,\beta}
    &=
    -M^2
    \sum_{\sigma=\alpha,\beta}
    \iint_{x\neq y}
    \left(
    \frac{x^2+y^2}{8}
    -\log|x-y|
    \right)
    d\mu_{\sigma,M}(x)d\mu_{\sigma,M}(y)
    \\
    &\quad
    +\sum_{\sigma=\alpha,\beta}
    \int
    \frac14 x^2\,d\mu_{\sigma,M}(x)
    \left(
    \operatorname{mass}(\mu_{\sigma,M})-\frac1M
    \right).
\end{align*}
Furthermore, we have
\begin{align*}      
   I_2^{\alpha,\beta}+&I_2^{\alpha}+I_2^{\beta}=
    -\sum_{\sigma=\alpha,\beta}\int \frac{1}{2}(\sqrt{2}+\Theta_M x)^2
    d\mu_{\sigma,M}(x)
    \nonumber\\
    &\quad
    +\iint
    \log(2\sqrt{2}+\Theta_M(x-y))^2
    d\mu_{\alpha,M}(x)d\mu_{\beta,M}(y)
    \nonumber\\
    &\quad=
    -\sum_{\sigma=\alpha,\beta}\int \frac{1}{2}(\sqrt{2}+\Theta_M x)^2
    d\mu_{\sigma,M}(x)
    \nonumber\\
    &\quad
    -\iint
    \left(
    \frac{x^2+y^2}{4}
    -
    \log(2\sqrt{2}+\Theta_M(x+y))^2
    \right)
    d\mu_{\alpha,M}(x)d\mu_{\beta,M}(y)
    \nonumber\\
    &\quad
    +\int\frac14x^2d\mu_{\alpha,M}(x)
    \operatorname{mass}(\mu_{\beta,M})
    +
    \int\frac14x^2d\mu_{\beta,M}(x)
    \operatorname{mass}(\mu_{\alpha,M}).
\end{align*}
Now, we deal with the terms $I_1^{\alpha},I_1^{\beta}$. By the definition of $\alpha_i,$ we have $\Theta_M\alpha_i+\sqrt{2}\ge 0$ and similarly $-\Theta_M\beta_j+\sqrt{2}\ge 0$. Moreover, observe that
\begin{equation*}
\left(\frac{\Theta_M}{\sqrt{2}}x\right)^2
+
2\frac{\Theta_M}{\sqrt{2}}x
-
2\log\left(1+\frac{\Theta_M}{\sqrt{2}}x\right)
=
f\left(\frac{\Theta_M}{\sqrt{2}}x\right),
\end{equation*}
where $f(x)=x^2+2x-2\log(1+x).$
One verifies that $f(x)\geq0$ for $x>-1$, and we have the pointwise convergence, $\Theta_M^{-2}
f\left(\frac{\Theta_M}{\sqrt{2}}x\right)
\longrightarrow x^2$ as $M\to\infty$. Let $R>0$, and denote 
$$
f_{\Theta_M,R}=\left[\Theta_M^{-2}f\left(\frac{\Theta_M}{\sqrt{2}}x\right)\right]\wedge R
$$ 
the truncated function, and one checks that $f_R$ is bounded continuous and $C\sqrt{R}$-Lipschitz with some constant $C>0$. Therefore, using the formula \eqref{eq:con GUE density alpha beta}, we obtain
\begin{align*}
&\quad \tilde{Q}( A_{\alpha,\beta,M}^{\varepsilon})\\
&\le
\sum_{k=0}^{M}
\int_{(A_{\alpha,\beta,M}^{\varepsilon})_k}
\exp\left[
-\Lambda_M\left(
\sum_{\sigma=\alpha,\beta}
\int
f_{\Theta_M,R}(x)
d\mu_{\sigma,M}(x)-\iint d\mu_{\alpha,M}(x)d\mu_{\beta,M}(y)\right)
\right]\\
&\quad\times \exp\left[
-M^2
(I_2^{\alpha}
+I_2^{\beta}
+I_2^{\alpha,\alpha}
+I_2^{\beta,\beta}
+I_2^{\alpha,\beta})
\right]
\prod_{\substack{1\leq i\leq k\\k+1\leq j\leq M}}d\alpha_id\beta_j.\\
&=
\exp\left[-\Lambda_M\left(
\sum_{\sigma=\alpha,\beta}
\int f_{\Theta_M,R}(x)
\,d\mu_\sigma(x) -\iint d\mu_\alpha(x)d\mu_\beta(y)+ O_R(\varepsilon)\right)\right]
\\
&\quad\times
\sum_{k=0}^{M}
\int_{(A_{\alpha,\beta,M}^{\varepsilon})_k}
\exp\Big[
-M^2
(
I_2^{\alpha}
+I_2^{\beta}
+I_2^{\alpha,\alpha}
+I_2^{\beta,\beta}
+I_2^{\alpha,\beta}
)
\Big]
\prod_{\substack{1\leq i\leq k\\k+1\leq j\leq M}}d\alpha_id\beta_j.
\end{align*}
where $A_{\alpha,\beta,M,k}^\varepsilon=(A_{\alpha,\beta,M}^\varepsilon)_k$ is defined as in \eqref{eq:con GUE density alpha beta}. The integral involving the terms of order $M^2$ is negligible at the large deviation scale $\Lambda_M=M^2\log\Theta_M^{-2},$ by the same argument used in the proof of the upper bound in \eqref{ineq:LDP sym block upperbound}. Moreover, since $f_{\Theta_M,R}\geq0$ and
$f_{\Theta_M,R}
\longrightarrow x^2\wedge R,$ Fatou's lemma yields
\begin{align*}
   \liminf_{M\to +\infty}\int f_{\Theta_M,R}(x)\,d\mu_{\sigma}(x)\ge \int x^2\wedge R\,d\mu_{\sigma}(x).
\end{align*}
Thus, we have for any $R>0$,
\begin{align*}
    \lim_{\varepsilon\to 0}\limsup_{M\to +\infty}\frac{1}{\Lambda_M}\log\,&\tilde{Q}_{M,N}( A_{\alpha,\beta,M}^{\varepsilon})\\
    &\le -\sum_{\sigma=\alpha,\beta}\int x^2\wedge R\,d\mu_\sigma(x)+\iint 1\,d\mu_\alpha(x)d\mu_\beta(y).
\end{align*}
Finally, pass the limit $R\to +\infty$, we obtain the desired inequality \eqref{ineq:LDP con GUE alpha beta limsup}. As a consequence, after absorbing the constant terms arising in the above discussion into the normalization constant $\tilde{Z}_{M,N}$, we obtain
\begin{equation*}
    \limsup_{M\to\infty}
    \frac{1}{\Lambda_M}
    \log \tilde{Z}_{M,N}
    \leq
    -\inf_{(\mu_\alpha,\mu_\beta)\in\mathcal{M}}
    I_{\alpha,\beta}(\mu_\alpha,\mu_\beta).
\end{equation*}

\medskip
\noindent\textbf{Lower bound.} We shall prove that, for every $(\mu_{\alpha},\mu_{\beta})\in\mathcal{M}$,
\begin{align}\label{ineq:LDP con GUE alpha beta liminf}
    \lim_{\varepsilon\to0}\liminf_{M\to\infty}
    \frac{1}{\Lambda_M}
    \log \tilde{Q}_{M,N}( A_{\alpha,\beta,M}^{\varepsilon})
    \geq
    -I_{\alpha,\beta}(\mu_\alpha,\mu_\beta).
\end{align}

As in the proof of the lower bound in the previous section, it suffices to establish \eqref{ineq:LDP con GUE alpha beta liminf} under the additional assumption that $\mu_{\alpha}$ and $\mu_{\beta}$ are compactly supported and admit continuous densities.

Assume moreover that $\operatorname{mass}(\mu_{\alpha})=p,\,
\operatorname{mass}(\mu_{\beta})=1-p,\,p\in[0,1].$
For $1\leq i\leq \lfloor Mp\rfloor$, define recursively
\begin{align*}
    x^{1,M}
    &=
    \inf\left\{
    x:
    \mu_{\alpha}((-\infty,x])
    \geq
    \frac{p}{M+1}
    \right\},
    \\
    x^{i+1,M}
    &=
    \inf\left\{
    x\geq x^{i,M}:
    \mu_{\alpha}((x^{i,M},x])
    \geq
    \frac{p}{M+1}
    \right\}.
\end{align*}
Similarly, for $\lfloor Mp\rfloor+1\leq j\leq M$, define
\begin{align*}
    x^{\lfloor Mp\rfloor+1,M}
    &=
    \inf\left\{
    x:
    \mu_{\beta}((-\infty,x])
    \geq
    \frac{1-p}{M+1}
    \right\},
    \\
    x^{j+1,M}
    &=
    \inf\left\{
    x\geq x^{j,M}:
    \mu_{\beta}((x^{j,M},x])
    \geq
    \frac{1-p}{M+1}
    \right\}.
\end{align*}
We then define, for $i=1,\ldots,\lfloor Mp\rfloor,\ j=\lfloor Mp\rfloor +1,\ldots,M$,
\[
A_{\alpha,\beta,M}^{\varepsilon,p}
=
\left\{
(\lambda_i)_{i=1}^M:
|\alpha_i-x^{i,M}|<\frac{\varepsilon}{2},
\;
|\beta_j+x^{j,M}|<\frac{\varepsilon}{2}
\right\}.
\]
By the same argument as in the proof of the lower bound in the previous section, for all sufficiently large $M$, $A_{\alpha,\beta,M}^{\varepsilon,p}
\subset A_{\alpha,\beta,M}^{\varepsilon}.$
Note that the points $x^{i,M}$, $i=1,\ldots,M$, are uniformly bounded. Hence, for every $\varepsilon>0$, there exists $M(\varepsilon)$ such that, for all $M\geq M(\varepsilon)$,
\begin{equation*}
    \log\left(1+\frac{\Theta_M x}{\sqrt{2}}\right)
    =
    \frac{\Theta_M x}{\sqrt{2}}
    -\frac{\Theta_M^2x^2}{4}
    +\theta\!\left(\frac{\Theta_M x}{\sqrt{2}}\right),
\end{equation*}
where $|\theta(y)|\leq |y|^3.$
Therefore, in $A_{\alpha,\beta,M}^{\varepsilon,p}$, 
\begin{align*}
    &\quad
    \prod_{i:\lambda_i\geq0}
    \left[
    (\sqrt{2}+\Theta_M\alpha_i)^{2(N-M)}
    \exp\left(
    -\frac{N-M}{2}
    (\sqrt{2}+\Theta_M\alpha_i)^2
    \right)
    \right]
    \\
    &=
    \exp\Bigg[
    -(N-M)
    \sum_{i=1}^{\lfloor Mp\rfloor}
    \left(
    \frac12(\sqrt2+\Theta_M\alpha_i)^2
    -
    2\log(\sqrt2+\Theta_M\alpha_i)
    \right)
    \Bigg]
    \\
    &=
    \exp\Bigg[
    -(N-M)
    \sum_{i=1}^{\lfloor Mp\rfloor}
    \left(
    1-\log2
    +\Theta_M^2\alpha_i^2
    -2\theta\!\left(
    \frac{\Theta_M\alpha_i}{\sqrt2}
    \right)
    \right)
    \Bigg]
    \\
    &\geq
    \exp\Bigg[
    -(N-M)\Theta_M^2
    \sum_{i=1}^{\lfloor Mp\rfloor}
    \left(
    \alpha_i^2
    +
    \frac{\Theta_M}{\sqrt2}
    |\alpha_i|^3
    \right)
    \Bigg]
    \\
    &\quad\times
    \exp\Big[
    -(N-M)\lfloor Mp\rfloor
    (1-\log2)
    \Big].
\end{align*}
Similarly,
\begin{align*}
    &\quad
    \prod_{i:\lambda_i<0}
    \left[
    (-\sqrt2+\Theta_M\beta_j)^{2(N-M)}
    \exp\left(
    -\frac{N-M}{2}
    (-\sqrt2+\Theta_M\beta_j)^2
    \right)
    \right]
    \\
    &\geq
    \exp\Bigg[
    -(N-M)\Theta_M^2
    \sum_{j=\lfloor Mp\rfloor+1}^{M}
    \left(
    \beta_j^2
    +
    \frac{\Theta_M}{\sqrt2}
    |\beta_j|^3
    \right)
    \Bigg]
    \\
    &\quad\times
    \exp\Big[
    -(N-M)
    (M-\lfloor Mp\rfloor)
    (1-\log2)
    \Big].
\end{align*}
Again, since
\begin{equation*}
    \exp\left[(N-M)\lfloor Mp \rfloor(1-\log 2)\right]
    \cdot
    \exp\left[(N-M)(M-\lfloor Mp \rfloor)(1-\log 2)\right]
\end{equation*}
depends only on $M$, it can be absorbed into the normalization constant $\tilde{Z}_{M,N}$.

Following a similar argument as in \eqref{ineq:LDP sym block lowerbound}, we obtain
\begin{align*}
\frac{1}{M}
\sum_{i=1}^{\lfloor Mp\rfloor}
\left(
\alpha_i^2
+
\frac{\Theta_M}{\sqrt{2}}
|\alpha_i|^3
\right)
&\sim
\frac{1}{M}
\sum_{i=1}^{\lfloor Mp\rfloor}
\left(
(x^{i,M})^2
+
\frac{\Theta_M}{\sqrt{2}}
|x^{i,M}|^3
\right)
\\
&\longrightarrow
\int x^2\,d\mu_{\alpha}(x),
\end{align*}
and the same convergence holds for the corresponding term involving $\beta_j$.

Moreover, on the set $A_{\alpha,\beta,M}^{\varepsilon,p}$, since $k=\lfloor Mp\rfloor$ we have
\[
k(M-k)\log\Theta_M^{-2}
\sim
p(1-p)\Lambda_M .
\]
Finally, for the terms
\begin{align*}
\prod_{1\leq i<j\leq \lfloor Mp\rfloor}
|\alpha_i-\alpha_j|^2
\cdot
\prod_{\lfloor Mp\rfloor\leq i<j\leq M}
|\beta_i-\beta_j|^2\cdot
\prod_{\substack{
1\leq i\leq \lfloor Mp\rfloor\\
\lfloor Mp\rfloor\leq j\leq M
}}
|2\sqrt{2}+\Theta_M(\alpha_i-\beta_j)|^2 ,
\end{align*}
appearing in \eqref{eq:con GUE density alpha beta}, the argument is identical to the one used for the lower bound in \eqref{ineq:LDP sym block lowerbound}. Hence, the integral of these terms over $A_{\alpha,\beta,M}^{\varepsilon,p}$ is of order $M^2$, which is negligible compared with the large deviation speed $\Lambda_M\gg M^2$. Therefore,
\begin{align*}
\lim_{\varepsilon\to0}
\liminf_{M\to\infty}
\frac{1}{\Lambda_M}
\log\tilde{Q}_{M,N}(A_{\alpha,\beta,M}^{\varepsilon})
&\geq
\lim_{\varepsilon\to0}
\liminf_{M\to\infty}
\frac{1}{\Lambda_M}
\log\hat{Q}_{M,N}(A_{\alpha,\beta,M}^{\varepsilon,p})
\\
&\geq
\iint d\mu_{\alpha}(x)d\mu_{\beta}(y)
-
\sum_{\sigma=\alpha,\beta}
\int x^2\,d\mu_{\sigma}(x).
\end{align*}

As a consequence, for every
$(\mu_{\alpha},\mu_{\beta})\in\mathcal{M}$,
\begin{equation*}
\liminf_{M\to\infty}
\frac{1}{\Lambda_M}
\log\tilde{Z}_{M,N}
\geq
\lim_{\varepsilon\to0}
\liminf_{M\to\infty}
\frac{1}{\Lambda_M}
\log\tilde{Q}_{M,N}(A_{\alpha,\beta,M}^{\varepsilon})
\geq
-I_{\alpha,\beta}(\mu_{\alpha},\mu_{\beta}),
\end{equation*}
which implies
\begin{align*}
\liminf_{M\to\infty}
\frac{1}{\Lambda_M}
\log\tilde{Z}_{M,N}
\geq
-\inf_{(\mu_{\alpha},\mu_{\beta})\in\mathcal{M}}
I_{\alpha,\beta}(\mu_{\alpha},\mu_{\beta}).
\end{align*}
\end{proof}

\medskip
\noindent\textbf{Conclusion.}
A simple observation is that the unique minimizer of the rate function
$I_{\alpha,\beta}$ over the set $\mathcal{M}$ is given by $(\mu_\alpha,\mu_\beta)
=
\left(\frac12\delta_0,\frac12\delta_0\right).$
Indeed, by the previous results, we have
\begin{align*}
   \inf_{(\mu_{\alpha},\mu_{\beta})\in\mathcal{M}}
   I_{\alpha,\beta}(\mu_{\alpha},\mu_{\beta})
   &=
   \inf_{p\in[0,1]}
   \inf_{\substack{\mathrm{mass}(\mu_\alpha)=p\\
   \mathrm{mass}(\mu_\beta)=1-p}}
   I_{\alpha,\beta}(\mu_\alpha,\mu_\beta)\\
   &=
   \inf_{p\in[0,1]}
   I_{\alpha,\beta}(p\delta_0,(1-p)\delta_0).
\end{align*}
The last quantity is minimized at $p=1/2$, which proves the claim. As a consequence, we obtain the following result.
\begin{corollary}
The pair
\[
(\mu_{\alpha,M},\mu_{\beta,M})
\stackrel{\mathcal{L}}{\longrightarrow}
\left(\frac12\delta_0,\frac12\delta_0\right)
\]
almost surely as $M\to+\infty$.
\end{corollary}

Consequently, we obtain the desired result.
\begin{corollary}
Almost surely, $\mathrm{mass}(\mu_{\alpha,M})
=
\mathrm{mass}(\mu_{\beta,M})
=1/2$ as $M\to\infty$. Combining this with Corollary \ref{cor:con GUE a.s.}, we obtain
\[
L_M
\stackrel{\mathcal{L}}{\longrightarrow}
\mu_{1/2}
=
\frac12\delta_{-\sqrt2}
+
\frac12\delta_{\sqrt2}
\]
almost surely as $M\to+\infty$.
\end{corollary}

\subsection{Some extensions}
\noindent\textbf{In the regime
\texorpdfstring{$M/N\to\alpha\in(0,1]$}
{M/N -> alpha in (0,1]}.} We now establish a large deviation principle for the empirical spectral measure $L_M$ associated with the density in \eqref{eq:con GUE density} in the regime $M/N\to\alpha\in(0,1]$. The proof follows essentially the same strategy as in \cite{arous1997large} and thus we omit it.

\begin{corollary}\label{cor:LDP con GUE alpha}
Assume that $M/N\to\alpha\in(0,1].$ Then the empirical spectral measure $L_M$ satisfies a large deviation principle with speed $NM$ and good rate function on $\mathcal{M}(\mathbb{R})$ given by
\begin{align*}
I_{\alpha}(\mu)
&=
\alpha\Sigma(\mu)
-\int
\left(
\frac{x^2}{2}
-(1-\alpha)\log x^2
\right)
d\mu(x)
-c_{\alpha},
\end{align*}
where
\[
c_{\alpha}
:=
\inf_{\mu\in\mathcal{M}(\mathbb{R})}
\left\{
\alpha\Sigma(\mu)
-\int
\left(
\frac{x^2}{2}
-(1-\alpha)\log x^2
\right)
d\mu(x)
\right\}.
\]
\end{corollary}

The rate function $I_{\alpha}$ provides an interpolation between the free entropy and the Boolean entropy relative to quadratic potential. Moreover, the minimizer of $I_{\alpha}$ admits an explicit density.

\begin{proposition}
There exists a unique probability measure $\mu_{\alpha}\in\mathcal{M}(\mathbb{R})$ such that $I_{\alpha}(\mu_{\alpha})=0.$ Its density is given by
\[
p_{\alpha}(x)
=
\frac{
\sqrt{(x^2-\gamma_1^2(\alpha))
(\gamma_2^2(\alpha)-x^2)}
}
{2\pi\alpha |x|}
\1_{[-\gamma_2(\alpha),-\gamma_1(\alpha)]
\cup
[\gamma_1(\alpha),\gamma_2(\alpha)]}(x),
\]
where $\gamma_1^2(\alpha)=
2-2\sqrt{2\alpha-\alpha^2},\ 
\gamma_2^2(\alpha)
=
2+2\sqrt{2\alpha-\alpha^2}.$
\end{proposition}

\begin{proof}
    The minimization problem is equivalent to solving the corresponding Euler--Lagrange equation for $I_\alpha$ (see e.g. \cite[Theorem 5.3.3]{hiai2000semicircle}). Namely, we look for a probability measure $\mu_\alpha\in\mathcal{M}(\mathbb{R})$ supported on $I_{a,b}:=[-\sqrt b,-\sqrt a]\cup[\sqrt a,\sqrt b]$ such that
    \begin{equation}\label{eq:euler lagrange minimizer}
       2\alpha\int\log|x-y|\,d\mu_\alpha(y)
       \begin{cases}
          =\frac{x^2}{2}-(\alpha-1)\log x^2+C,
          & x\in I_{a,b},\\
          \leq \frac{x^2}{2}-(\alpha-1)\log x^2+C,
          & \text{otherwise},
       \end{cases}
    \end{equation}
    where $C$ is a constant. Observe that
\[
\int\log|x-y|\,d\mu_\alpha(y)
=
\lim_{\varepsilon\to0^+}
\Re
\int\log(x+i\varepsilon-y)\,d\mu_\alpha(y)
=
\lim_{\varepsilon\to0^+}\Re F(x+i\varepsilon),
\]
and
\[
F'(z)
=
\int\frac{1}{z-y}\,d\mu_\alpha(y)
=
G_{\mu_\alpha}(z),
\]
where $G_{\mu_\alpha}$ is the Cauchy transform of $\mu_\alpha$, analytic on $\mathbb{C}^+=\{z\in\mathbb{C}:\Im z>0\}.$ Therefore, a necessary condition for $\mu_\alpha$ to satisfy the Euler--Lagrange equation is
\begin{equation}\label{eq:necessary condition}
\lim_{\varepsilon\to0^+}
\Re G_{\mu_\alpha}(x+i\varepsilon)
=
\frac{1}{2\alpha}
\left(
x-(\alpha-1)\frac{2}{x}
\right),
\qquad x\in I_{a,b}.
\end{equation}
The left-hand side is precisely the Hilbert transform of $\mu_\alpha$ up to a normalization factor.

We recall the following standard fact \cite[Chapter 3, Theorem 10]{mingo2017free}: if $H:\mathbb{C}^+\to\mathbb{C}^-$
is analytic and $\limsup_{y\to+\infty}y|H(iy)|=c<\infty,$ then there exists a unique positive Borel measure $\nu$ on $\mathbb{R}$ such that
\[
H(z)=\int\frac{1}{z-x}\,d\nu(x),
\qquad
\nu(\mathbb{R})=c.
\]

Motivated by the computation of the Marchenko--Pastur distribution $\nu_\gamma$ \eqref{eq:MP density} as the minimizer of the rate function $J_\gamma^+$ \eqref{eq:rate function wishart} (see \cite[Section 5.5]{hiai2000semicircle}), together with the necessary condition \eqref{eq:necessary condition}, we consider the family of analytic functions
\[
H_{a,b}(z)
=
\frac{z^2+2(1-\alpha)-\sqrt{(z^2-a)(z^2-b)}}{2\alpha z},
\]
parametrized by $b>a\ge0$.

We choose the branch such that for $x\in[\sqrt b,\infty)$,
\[
\sqrt{(x^2-a)(x^2-b)}\ge0.
\]
We first determine $a,b$ such that $H_{a,b}$ is the Cauchy transform of a probability measure $\mu_\alpha$. For parameters satisfying $ab=4(\alpha-1)^2,$
we have $H_{a,b}:\mathbb{C}^+\to\mathbb{C}^-.$
Moreover,
\[
\lim_{y\to+\infty}y|H_{a,b}(iy)|
=
\frac{4(\alpha-1)+a+b}{4\alpha}.
\]
Since we require the corresponding measure to be a probability measure, we impose
\[
\frac{4(\alpha-1)+a+b}{4\alpha}=1.
\]
Together with the previous relation, this gives
\[
\begin{cases}
a+b=4,\\
ab=4(\alpha-1)^2,
\end{cases}
\]
and hence $a=2-2\sqrt{2\alpha-\alpha^2},\ b=2+2\sqrt{2\alpha-\alpha^2}.$ The corresponding density is obtained from the Stieltjes inversion formula:
\[
p_\alpha(x)
=
\lim_{\varepsilon\to0^+}
\frac1\pi
\Im H_{a,b}(x+i\varepsilon),
\]
which yields the desired density $p_\alpha$.

Finally, it remains to check \eqref{eq:euler lagrange minimizer}, setting $F'(z)=H_{a,b}(z),$ we obtain
\[
\frac1\pi
\lim_{\varepsilon\to0^+}
\Re H_{a,b}(x+i\varepsilon)
=
\frac1{2\alpha\pi}
\left(
x+(1-\alpha)\frac2x
\right),
\]
for $x\in I_{a,b}$. Hence, in the sense of distributions,
\[
F'(x)
=
\frac1{2\alpha}
\left(
x-(\alpha-1)\frac2x
\right),
\]
and therefore
\[
F(x)
=
\frac1{2\alpha}
\left(
\frac{x^2}{2}
-(\alpha-1)\log x^2
\right)
+C
\]
on $I_{a,b}$. For $x\notin I_{a,b}$, the function $F$ is differentiable and
\[
F'(x)
=
\int\frac{p_\alpha(y)}{x-y}dy
\]
is given explicitly by the expression for $H_{a,b}$. Without loss of generality, let $x>\sqrt b$, we have
\[
F'(x)
<
\frac{x^2+2(1-\alpha)}{2\alpha x}.
\]
Since $F$ is continuous at $\sqrt b$, we conclude that
\[
F(x)
<
\frac1{2\alpha}
\left(
\frac{x^2}{2}
-(\alpha-1)\log x^2
\right)+C.
\]
which completes the proof.
\end{proof}

\begin{remark}
Indeed, the uniqueness of the minimizer of $I_\alpha$ in Corollary \ref{cor:LDP con GUE alpha} is mainly a consequence of the free entropy term $\Sigma(\mu)$ appearing in the rate function $I_1(\mu)$. In particular, $\Sigma(\mu)$ is strictly concave on the set of probability measures supported on any fixed compact subset of $\mathbb{C}$ (see e.g. \cite[Proposition 5.3.2]{hiai2000semicircle}). This quantity originates from the Vandermonde term
\[
\prod_{1\leq i<j\leq M}|\lambda_i-\lambda_j|^2,
\]
which encodes the interaction between eigenvalues. The presence of this term is also essential in the proof of Theorem \ref{thm:LDP rescaled con GUE}.
\end{remark}

\begin{remark}
When $\alpha=1$, one can verify that
\[
p_1(x)=\frac{1}{2\pi}\sqrt{4-x^2}\,\1_{[-2,2]}(x),
\]
which recovers the semicircle law as the unique minimizer of $I_1$. On the other hand, as $\alpha\to0$, one can check that
\[
p_\alpha(x)\longrightarrow
\mu_{1/2}
=
\frac12\delta_{-\sqrt2}
+
\frac12\delta_{\sqrt2}
\]
in distribution.
\end{remark}

\noindent\textbf{Models with general potential and beta.} We can generalize the density in \eqref{eq:con GUE density} as follows (up to a normalized constant)
\begin{equation}\label{eq:general con GUE density}
    \prod_{i=1}^M |\lambda_i|^{\gamma_N}
    \prod_{i<j}|\lambda_i-\lambda_j|^{2\beta_N}
    \prod_{i=1}^M e^{-N V(\lambda_i)}d\lambda_i,
\end{equation}
where $\beta_N>0,\ \gamma_N\geq 0,\ M=M_N$ depend on $N$. We assume that $V:\mathbb{R}\to\mathbb{R}$ is a continuous function satisfying
\[
\lim_{x\to\infty} xe^{-\varepsilon V(x)}=0
\]
for every $\varepsilon>0$. Let $M\to \infty,\ M/N\to0$ as $N\to +\infty$. We still denote the empirical spectral measure by
\[
L_M=\frac{1}{M}\sum_{i=1}^M\delta_{\lambda_i}.
\]
By adapting the proof of Theorem \ref{thm:LDP con GUE} and Corollary \ref{cor:LDP con GUE alpha}, we obtain the following result.

\begin{corollary}\label{cor:LDP con beta ensemble}
Assume further that $\gamma_N/N\to\gamma\ge0$ and if $\beta_N=o(N/M)$
as $N\to\infty$. Then $L_M$ satisfies a large deviation principle with speed $MN$ and good rate function
\[
I_{\gamma,V}(\mu) =
\int (V(x)-\gamma\log|x|)\,d\mu(x)
-c_{\gamma,V},
\]
where
$$
c_{\gamma,V} = \inf_{\mu\in\mathcal{M}(\mathbb{R})}\left\{
\int (V(x)-\gamma\log|x|)\,d\mu(x)\right\}.
$$
The minimizer of $I_{\gamma,V}$ is not necessarily unique.

Else if $\beta_N\cdot M/N\to \beta>0$ as $N\to\infty$, then $L_M$ satisfies a large deviation principle with speed $MN$ and good rate function
\[
I_{\gamma,\beta,V}(\mu) = -\beta\Sigma(\mu)+
\int (V(x)-\gamma\log|x|)\,d\mu(x)
-c_{\gamma,\beta,V},
\]
where
$$
c_{\gamma,V} = \inf_{\mu\in\mathcal{M}(\mathbb{R})}\left\{-\beta\Sigma(\mu)+
\int (V(x)-\gamma\log|x|)\,d\mu(x)\right\}.
$$
The minimizer of $I_{\gamma,\beta,V}$ is unique.
\end{corollary}

A similar phenomenon appears in Girko's theorem for the complex Ginibre ensemble. More precisely, if $X(N)$ is a random matrix whose entries are i.i.d. complex random variables with mean zero and variance $1$, then the empirical spectral measure of $X_N=N^{-1/2}X(N)$ converges almost surely weakly to the circular law, namely the uniform probability measure on the unit disk in $\mathbb{C}$. 

In the Gaussian case, and for simplicity restricting to real entries, the joint density of the eigenvalues is given by
\begin{equation*}
    dU_N(\zeta_1,\ldots,\zeta_N)
    =
    \frac{1}{C_N}
    \prod_{i<j}|\zeta_i-\zeta_j|^2
    \prod_{i=1}^N e^{-N|\zeta_i|^2}d\zeta_i .
\end{equation*}
This expression also leads to a large deviation principle for the empirical spectral measure at speed $N^2$ with good rate function
\begin{equation*}
    I^c(\mu)
    :=
    -\Sigma(\mu)
    +
    \int |\zeta|^2\,d\mu(\zeta)
    -
    b_0^c,
    \qquad
    \mu\in\mathcal{M}(\mathbb{C}),
\end{equation*}
where
\[
b_0^c
:=
\inf_{\nu\in\mathcal{M}(\mathbb{C})}
\left\{
-\Sigma(\nu)
+
\int |\zeta|^2\,d\nu(\zeta)
\right\}.
\]
Moreover, $I^c$ admits a unique minimizer, which is precisely the circular law. We refer the reader to \cite[Sec.~5.4]{hiai2000semicircle} for full treatment. Conditioning $N-M$ eigenvalues to be zero, the remaining $M$ eigenvalues have a density of the same type as \eqref{eq:con GUE density},
\begin{equation}\label{eq:con Ginibre density}
    dU_{M,N}(\zeta_1,\ldots,\zeta_M)
    =
    \frac{1}{C_{M,N}}
    \prod_{i=1}^M|\zeta_i|^{2(N-M)}
    \prod_{ i<j\leq }|\zeta_i-\zeta_j|^2
    \prod_{i=1}^M e^{-N|\zeta_i|^2}d\zeta_i .
\end{equation}
Consequently, one obtains a corresponding large deviation principle with speed $NM$.

\begin{corollary}\label{cor:LDP con Ginibre}
In the regime $M/N\to \alpha\in [0,1]$ as $N\to \infty$. Under the law $U_{M,N}$, the empirical spectral measure on $\mathbb{C}$ satisfies a large deviation principle with speed $NM$ and good rate function
\[
I_\alpha^c(\mu) = -\alpha\iint\log|\zeta-\zeta'|\,d\mu(\zeta)d\mu(\zeta')+
\int
\left(
|\zeta|^2-(1-\alpha)\log|\zeta|^2
\right)
d\mu(\zeta)
-c_\alpha^c,
\]
where $c_\alpha^c$ denotes the infimum of the above functional over all $\mu\in\mathcal{M}(\mathbb{C})$.

\end{corollary}

\begin{remark}
When $\alpha=0$, the minimum of $I_\alpha^c$ is attained by any probability measure supported on the unit circle $S^1$. By symmetry, one expects the empirical spectral measure $L_M$ to converge to the uniform distribution on $S^1$. For $\alpha>0$, the minimizer is unique, and its support is expected to be contained in an annulus.
\end{remark}

In the regime $M/N\to 0$, a similar extension of the large deviation principle holds for the more general density proportional to
\[
\prod_{i=1}^M |\zeta_i|^{\gamma_N}
\prod_{i<j} |\zeta_i-\zeta_j|^{2\beta_N}
\prod_{i=1}^M e^{-N V(\zeta_i)}\, d\zeta_i,
\]
where $\gamma_N$ and $\beta_N$ are as in \eqref{eq:general con GUE density}, and $V:\mathbb{C}\to\mathbb{R}$
is continuous and satisfies, for every $\varepsilon>0$,
\[
\lim_{|\zeta|\to\infty}
|\zeta|e^{-\varepsilon V(\zeta)}
=0.
\]

Assume that $\gamma_N/N\to \gamma\ge0$ and $\beta_N M/N\to \beta\ge0$. Then the corresponding empirical spectral measure satisfies a large deviation principle with speed $NM$ and good rate function
\[
I_{\gamma,\beta,V}^c(\mu)
=
-\beta\iint \log|\zeta-\zeta'|\, d\mu(\zeta)\, d\mu(\zeta')
+
\int
\bigl(
V(\zeta)-\gamma\log|\zeta|^2
\bigr)\, d\mu(\zeta)
-
c_{\gamma,\beta,V}^c,
\]
where $c_{\gamma,\beta,V}^c$ denotes the infimum of the above functional over all $\mu\in\mathcal{M}(\mathbb{C})$.

\section{A Proof of Theorem \ref{thm:max and mono}}
In this section, we prove the maximality and monotonicity of the Boolean entropy defined by
\[
\Gamma(\mu)=\int \log |x|\,d\mu(x),\qquad \mu\in\mathcal{M}(\mathbb{R}),
\]
in analogy with the classical and free settings. Denote $\mathcal{P}^2$ is the set of probability measures with finite second moments and $\mathcal{P}_0^2$ stands for the space of probability measures with mean $0$ and finite variance. Furthermore, we denote by $\mathcal{P}_{0,c}^2$ the space of compactly supported probability measures with mean $0$ and finite variance. The maximality in Theorem \ref{thm:max and mono} is just a consequence Jensen's inequality.

For the monotonicity, to begin with, we introduce the following notation. Let $\{a_i\}_{i\ge 1}\subset(\mathcal{A},\varphi)$ be a sequence of identically distributed and Boolean independent random variables with a common distribution $\mu\in\mathcal{M}(\mathbb{R})$. We denote by
\[
\underbrace{\mu \uplus \cdots \uplus \mu}_{n \text{ times}}
\]
the distribution of the sum $a_1 + \cdots + a_n$. Moreover, for $\lambda > 0$, we denote by $D_\lambda(\mu)$ the distribution of the rescaled random variable $\lambda a_1$.  To prove the second part of Theorem \ref{thm:max and mono}, we first show the following.

\begin{proposition}\label{prop:monotonicity of entropy}
	For any $\mu\in \mathcal{P}_{0,c}^2$ such that $\Gamma(\mu)>-\infty$, and for any integers $n\ge m\ge 1$, we have
	\begin{equation}\label{ineq:mono entropy}
		\gamma_{\mu}(n):=\Gamma(\underbrace {D_{\frac{1}{\sqrt{n}}}\mu \uplus \cdots \uplus D_{\frac{1}{\sqrt{n}}}\mu}_{n\ \mathrm{times}})\ge \gamma_\mu(m).
	\end{equation}
	Moreover, $\{\gamma_{\mu}(n)\}_{n\ge 1}$ is a constant sequence if and only if $\mu= D_\lambda (\mathrm{b})$ for some $\lambda\ne 0$.
\end{proposition}
\begin{proof}
	Denote $\mu_n=\underbrace {D_{\frac{1}{\sqrt{n}}}\mu \uplus \cdots \uplus D_{\frac{1}{\sqrt{n}}}\mu}_{n\ \mathrm{times}}$ and recall that
\[
G_{\mu}(z)=\int_{\mathbb{R}} \frac{1}{z-x}\,d\mu(x),\qquad z\in\mathbb{C}_+
\]
is the Cauchy transform of a probability measure $\mu$ on $\mathbb{R}$, and the self-energy ($K$-transform) of $\mu$ is defined as
\[
K_\mu(z)=z-\frac{1}{G_{\mu}(z)}.
\]
By \cite[Proposition 2.1]{speicher1996boolean}, $K_{\mu}(z)$ is additive with respect to Boolean convolution, i.e.,
\begin{equation*}
K_{\mu_1 \uplus \mu_2}(z)=K_{\mu_1}(z)+K_{\mu_2}(z), \qquad \forall z\in\mathbb{C}^+.
\end{equation*}
Using this additivity and the scaling property $G_{D_\lambda\mu}(z)=\lambda^{-1}G_\mu(\lambda^{-1}z)$ for any $\lambda> 0$, we obtain
\begin{equation*}
G_{\mu_n}(z)
=
\frac{1}{z-\sqrt{n}\,K_{\mu}(\sqrt{n}z)}.
\end{equation*}
In the above expression, note that the parameter $n$ on the right-hand side can be formally replaced by any positive real number $t$, thereby defining the Cauchy transforms of a family of probability measures $(\mu_t)_{t\ge 1}$. We make this statement precise in the following lemma and postpone its proof.

\begin{lemma}\label{lem:extension}
Given a probability measure $\mu\in\mathcal{P}_{0,c}^2$, then there exists a family of probability measures $(\mu_t)_{t\geq 1}\subset \mathcal{P}_{0,c}^2$ extending the sequence $\{\mu_n\}_{n\in\mathbb{N}}$, defined by
\begin{equation}\label{eq:boolean CLT cauchy flow}
G_{\mu_t}(z) := G_t(z)
= \frac{1}{z-\sqrt{t}K_{\mu}(\sqrt{t}z)},
\end{equation}
and satisfies
\begin{enumerate}[leftmargin=20pt,label=(\alph*)]
    \item The family $(\mu_t)_{t\geq 1}$ forms a multiplicative semigroup, i.e., $(\mu_t)_s=\mu_{ts}$ for any $t,s\ge 1$.
    
    \item If $\mu \neq D_\lambda(\mathrm{b})$ for some $\lambda\ne 0$, then $\mu_t \neq D_\lambda(\mathrm{b})$ for all $t\ge 1$. Additionally, if $\mu(\{0\})=0,$ then for any $t\ge 1$, $\mu_t(\{0\})=0$.
\end{enumerate}
\end{lemma}
    We now aim to show that $\gamma_{\mu}(t):=\Gamma(\mu_t)$ is a non-decreasing function of $t$. It suffices to prove that $\gamma_{\mu}'(t)\ge 0$ for all $t\ge 1$, with equality if and only if $\mu=D_\lambda(\mathrm{b})$ for some $\lambda\ne 0$. Due to the singularity of $\log |\cdot|$ at the origin, we introduce a regularized logarithmic integral with parameter $\delta>0$, which makes sense to compute the derivative along $\mu_t$, and then pass to the limit as $\delta\to 0$ to obtain an expression for $\gamma_\mu'(t)$.

To be more specific, let $z=x+i\varepsilon$ with $\varepsilon>0$, and denote
$u_t(z)=\pi^{-1}\Re G_t(z),\  v_t(z)=-\pi^{-1}\Im G_t(z)$. By the definition of $\mu_t$ \eqref{eq:boolean CLT cauchy flow}, we have
\[
-\frac{\partial}{\partial t}\Big|_{t=1}v_t(z)
=\pi(u_1(z)^2-v_1(z)^2)\varepsilon+\frac{1}{2} u_1'(z)\varepsilon
-2\pi u_1(z)v_1(z)x+\frac{1}{2} v_1(z)-\frac{1}{2} v_1'(z)x.
\]
In addition, since $\mu$ is compactly supported, one has the bound with some constant $C_\varepsilon>0$,
\[
\left| \frac{\partial}{\partial t}\Big|_{t=1}v_t(z) \right|
\leq \frac{C_\varepsilon}{1+|x|^2}.
\]
For any $\delta\ge 0$, define $f_\delta(x)=\log |x+i\delta|,\ g_\delta(x)=x\log |x+i\delta|$ and 
$$
\Gamma_\delta(\mu)= \int f_\delta(x)\,d\mu(x).
$$
Note that $v_t(x+i\varepsilon)=P_\varepsilon\ast \mu_t(x)$, where $P_\varepsilon$ denotes the Poisson kernel. Hence,
\begin{equation}\label{eq: Poisson transform on entropy}
\Gamma(\mu_{t,\varepsilon})= \int (P_\varepsilon\ast \log)(x)\,d\mu_t(x)
=\int \log|x+i\varepsilon|\,d\mu_t(x)=\Gamma_\varepsilon(\mu_t),
\end{equation}
where $\mu_{t,\varepsilon}$ is the probability measure with density $v(\cdot+i\varepsilon)$.
With these observations, one computes,
\begin{align*}
\frac{d}{dt}\Big|_{t=1}\Gamma_{\delta+\varepsilon}(\mu_t)
&= \int \frac{\partial}{\partial t}\Big|_{t=1}v_t(x+i\varepsilon)\cdot \log |x+i\delta|\,dx \\
&= \int \left[2\pi u_1(z)v_1(z) +\frac{1}{2} v_1'(z)\right]g_\delta(x)\,dx - \int \frac{1}{2} v_1(z)f_\delta(x)\,dx \\
&\quad + \varepsilon \int \left[ \pi(u_1(z)^2-v_1(z)^2)+\frac{1}{2}u_1'(z)\right]f_\delta(x)\,dx \\
&\stackrel{\mathrm{def}}{=} R_{\delta,\varepsilon}(\mu_1).
\end{align*}
Thanks to the multiplicative semigroup property, plus $\mu_t$ is compactly supported for any $t\ge 1$ (by Lemma \ref{lem:extension}), following the same strategy, we have $\frac{d}{dt}\Gamma_{\delta+\varepsilon}(\mu_t) = \frac{1}{t}R_{\delta,\varepsilon}(\mu_t)$, which implies
\begin{equation}\label{eq:integral clt entropy delta varepsilon}
    \Gamma_{\delta+\varepsilon}(\mu_t) - \Gamma_{\delta+\varepsilon}(\mu)
= \int_1^t \frac{1}{s}R_{\delta,\varepsilon}(\mu_s)\,ds.
\end{equation}

We now claim a technical lemma describing behavior of $R_{\delta,\varepsilon}(\mu)$ as $\varepsilon\to 0$ and $\delta\to 0$ for any compactly supported probability measure $\mu$. The proof will be postponed.
\begin{lemma}\label{lem:tech lemma R_delta,varepsilon}
    For any compactly supported probability measure $\mu$, first for fixed $\delta>0$, we have
    \begin{equation}\label{eq:limit R_delta,varepsilon varepsilon to 0}
        \lim_{\varepsilon\downarrow0} R_{\delta,\varepsilon}(\mu) = \frac{1}{2}\iint\left[ \frac{x+y}{x-y}(f_\delta(x)-f_\delta(y))-\frac{x^2}{x^2+\delta^2} \right]\,d\mu(x)d\mu(y)
    \end{equation}
    uniformly in $\mu$. Denote
    \begin{equation*}
        R_\delta(\mu) = \iint \frac{x+y}{x-y}(f_\delta(x)-f_\delta(y))\,d\mu(x)d\mu(y),
    \end{equation*}
    we have $R_\delta(\mu)\ge 0$ for all $\delta\ge 0$ and as $\delta\downarrow 0^+$, the following monotone convergence holds
    \begin{equation}\label{eq:limit R_delta delta to 0}
         R_{\delta}(\mu) \nearrow R_0(\mu).
    \end{equation}
\end{lemma}

With these preparations, let first $\varepsilon\to 0$ in \eqref{eq:integral clt entropy delta varepsilon}, since $f_{\delta+\varepsilon}\le f_{\delta+1}$ for sufficiently small $\varepsilon$ and $\mu_t,\mu$ are compactly supported, the left-hand side converges to $\Gamma_\delta(\mu_t)-\Gamma_\delta(\mu)$ simply by the monotone convergence theorem. According to Lemma \ref{lem:extension}, each probability measure $\mu_s$ is compactly supported for any $s\ge 1$, thus applying \eqref{eq:limit R_delta,varepsilon varepsilon to 0}, together with the dominated convergence theorem on the right-hand side of \eqref{eq:integral clt entropy delta varepsilon}, yield
\begin{align*}
    \Gamma_\delta&(\mu_t)-\Gamma_\delta(\mu)= \int_1^t \frac{1}{2s}\left(R_\delta(\mu_s)-\int \frac{x^2}{x^2+\delta^2}\,d\mu_s(x)\right)\,ds.
\end{align*}
We now let $\delta \to 0$. Similarly the left-hand side converges to $\Gamma(\mu_t)-\Gamma(\mu)$. Applying Lemma~\ref{lem:tech lemma R_delta,varepsilon} to $\mu_s$, we obtain the monotone convergence of $R_\delta(\mu_s)$ as $\delta\to0$ (see \eqref{eq:limit R_delta delta to 0}) for every $s\geq 1$. Combining this with the fact that
$\frac{x^2}{x^2+\delta^2}\nearrow \1_{x\ne 0}$ as $\delta\to0$ and $\mu_s(\{0\})=0$ for any $s\ge 1$ (Lemma \ref{lem:extension}), we deduce that
\begin{equation}\label{eq:integral clt entropy}
    \Gamma(\mu_t)-\Gamma(\mu) = \int_1^t \frac{1}{2s}\left(R_0(\mu_s)-1\right)ds.
\end{equation}
Recall that $f_0(x)=\log |x|$, we have
\begin{align}\label{eq:derivative of entropy along clt}
    R_0(\mu_s)-1 =\iint \frac{x/y+1}{x/y-1}\log |x/y|\,d\mu_s(x)d\mu_s(y)-1.
\end{align}
Let $X,Y$ be classically i.i.d.\ random variables with law $\mu_s$, and define $Z=X/Y$. As $Z$ is well defined $\mu_s$-almost everywhere, the double integral in the equation \eqref{eq:derivative of entropy along clt} is equal to $\mathbb{E}[\ell(Z)]$, where $\ell(x)=\frac{x+1}{x-1}\log |x|.$ It can be checked that the function $\ell(x)$ is non-negative and has two local minima: $0$ and $2$, attained at $x=-1$ and $x=1$, respectively (see Figure \ref{fig:ell}).
\begin{figure}[ht]
    \centering
    \includegraphics[width=0.3\textwidth, height=4cm]{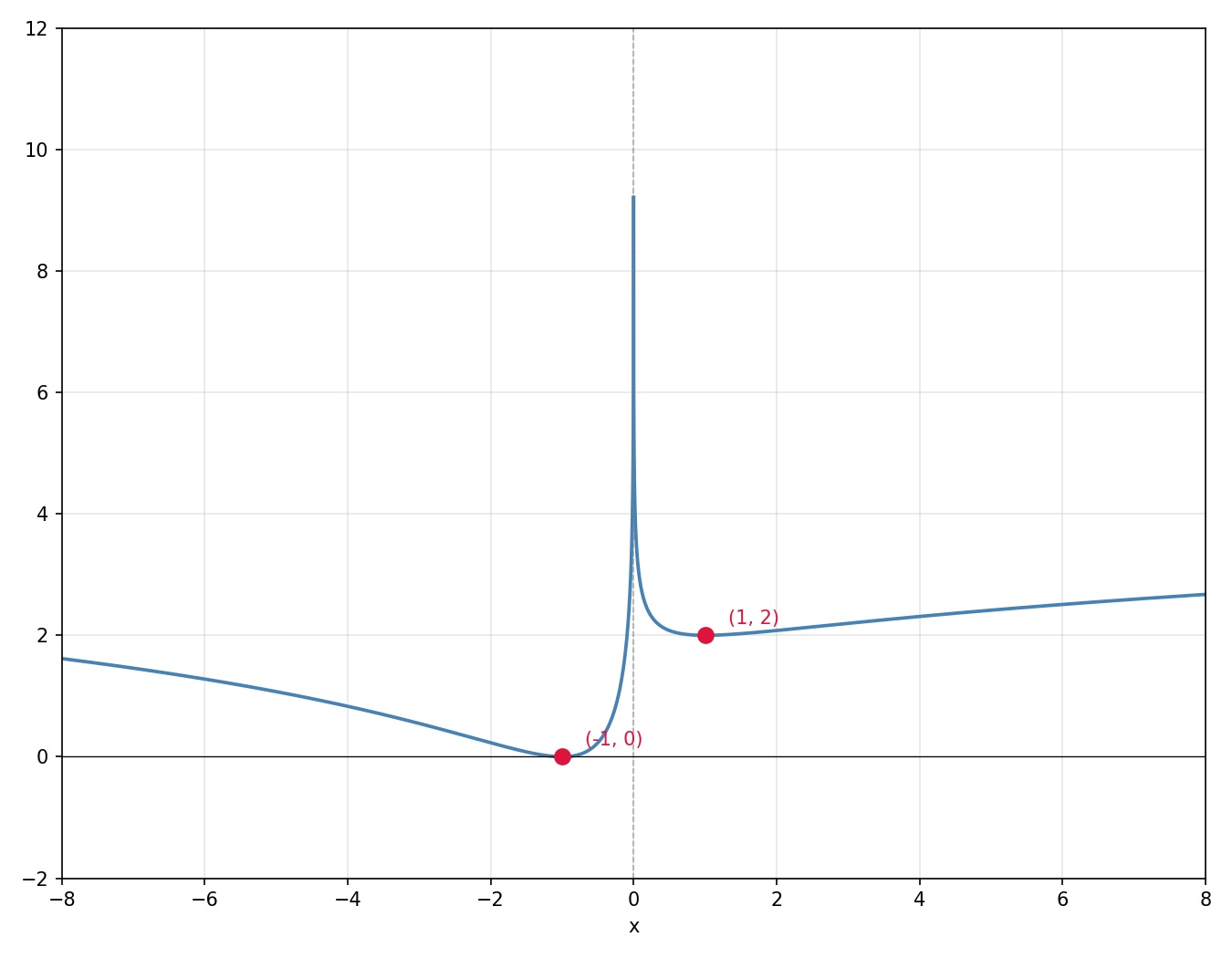}
    \caption{$y=\ell(x)$ with two local minima attained at $-1$ and $1$.}
    \label{fig:ell}
\end{figure}
Hence, let $\mathbb{P}(X\geq 0)=\theta\in[0,1]$, by the Markov inequality, we have
\begin{align*}
\mathbb{E}[\ell(Z)]
\geq \mathbb{E}[\ell(Z)\1_{\{Z\geq 0\}}]
&\geq 2\mathbb{P}(Z\geq 0)\\
&=2[\mathbb{P}(X\geq 0, Y\geq 0)+\mathbb{P}(X<0,Y<0)]\\
&=2[\theta^2+(1-\theta)^2]\geq 1.
\end{align*}
Moreover, the equality holds if and only if $\theta=1/2$ and $\mathbb{P}(Z=\pm1)=1$. By the assumption $\mu_s\in\mathcal{P}_{0,c}^2$, we deduce that $\mu_s=D_\lambda(\mathrm{b})$ for some $\lambda\ne 0$. Substituting this into \eqref{eq:derivative of entropy along clt} and using
\eqref{eq:integral clt entropy}, together with part (b) of
Lemma~\ref{lem:extension}, we obtain
$\gamma_\mu(t)\geq \gamma_\mu(s)$ for any $t\ge s\ge 1$, with equality if and only if
$\mu=D_\lambda(\mathrm{b})$ for some $\lambda\ne 0$. This completes the proof for monotonicity.
\end{proof}
\begin{remark}
    If $\Gamma(\mu)>-\infty$, by the monotonicity we know that $\Gamma(\mu_s)>-\infty$ for all $s\ge 1$, and since $\mu_s$ is compactly supported, we have $|\Gamma(\mu_s)|<\infty$. In particular, the random variable $Z$ is log integrable, i.e. $\mathbb{E}[|\log |Z||]<+\infty$, in view of \eqref{eq:derivative of entropy along clt}, we deduce that $R_0(\mu_s)<+\infty$ so that $\gamma'(s)=\frac{1}{s}(R_0(s)-\frac{1}{2}).$
\end{remark}

To complete the argument, we first prove Lemma \ref{lem:extension}.
\begin{proof}
    By \cite[Proposition 3.2]{speicher1996boolean}, for any $\mu\in\mathcal{P}_{0}^2$ with variance $\sigma^2>0$, there exists a positive measure $\rho$ on $\mathbb{R}$ with total mass $\sigma^2$ such that
\[
K_{\mu}(z)=\int \frac{1}{z-x}\,d\rho(x),
\]
and vice versa. By the scaling property of the $K$-transform, we have
\[
K_{\mu_t}(z)=\sqrt{t}\,K_{\mu}(\sqrt{t}z)
= \sqrt{t}\int \frac{1}{\sqrt{t}z-x}\,d\rho(x)
= \int \frac{1}{z-x/\sqrt{t}}\,d\rho(x).
\]
Let $\rho_t$ be the push-forward measure of $\rho$ with respect to the map $x\mapsto x/\sqrt{t}$. Then
\[
K_{\mu_t}(z)=\int \frac{1}{z-x}\,d\rho_t(x),
\]
which implies that $\mu_t \in \mathcal{P}_0^2$ with the same variance $\sigma^2$ for all $t\ge 1$. And it is clear that if $K_\mu(z)\ne K_{D_{\sigma^2}(\mathrm{b})}(z)=\sigma^2/z$, then for any $t\ge 1,$ $K_{\mu_t}(z)\ne \sigma^2/z$. Moreover, recall that $0$ is not an atom of $\mu$ if and only if
\[
\lim_{\varepsilon\to 0} i\varepsilon G_\mu(i\varepsilon)
=
\mu(\{0\})
=
0.
\]
Using the Cauchy transform of $\mu_t$ \eqref{eq:boolean CLT cauchy flow}, we have
\begin{align*}
i\varepsilon G_{\mu_t}(i\varepsilon)
=
\frac{i\varepsilon}
{\dfrac{\sqrt{t}}{G_\mu(\sqrt{t}i\varepsilon)}-(t-1)i\varepsilon}=
\frac{1}
{\dfrac{t}{\sqrt{t}i\varepsilon\,G_\mu(\sqrt{t}i\varepsilon)}-(t-1)}.
\end{align*}
Taking the limit as $\varepsilon\to 0$, we obtain
\[
\lim_{\varepsilon\to 0} i\varepsilon G_{\mu_t}(i\varepsilon)=0.
\]
This proves \textup{(b)}.

For the multiplicative semigroup property, note that for any $t \ge s \ge 1$, we have
\[
K_{\mu_t}(z)
= \sqrt{t}\,K_{\mu}(\sqrt{t}z)
= \sqrt{\frac{t}{s}}\,K_{\mu_s}\!\left(\sqrt{\frac{t}{s}}\,z\right),
\]
which means that $(\mu_s)_{t/s} = \mu_t.$ This proves (a).

It remains to show $\mu_t$ is compactly supported. Without loss of generality, assume that $\mu\in\mathcal P_{0,c}^2$ has unit variance ($\sigma^2=1$), by \cite[Proposition 2.1]{speicher1996boolean},
\begin{equation}\label{eq:K-transform generating function}
    K_\mu(z)
 =\sum_{k\ge1}\frac{r_n(\mu)}{z^{\,n-1}}
 =\frac1z+\sum_{k\ge3}\frac{r_n(\mu)}{z^{\,n-1}}
\end{equation}
for sufficiently large $|z|$, where $\{r_k(\mu)\}_{n\ge1}$ denotes the sequence of Boolean cumulants of $\mu$. Moreover, moments $\{m_k(\mu)\}_{k\ge 1}$ and Boolean cumulants satisfy the recursive relation
\[
m_k(\mu)
=
\sum_{p=1}^k
\sum_{\substack{i(1),\ldots,i(p)\ge1\\
i(1)+\cdots+i(p)=k}}
r_{i(1)}(\mu)\cdots r_{i(p)}(\mu),
\]
which can be inverted as
\[
r_k(\mu)
=
\sum_{p=1}^k
(-1)^{p+1}
\sum_{\substack{i(1),\ldots,i(p)\ge1\\
i(1)+\cdots+i(p)=k}}
m_{i(1)}(\mu)\cdots m_{i(p)}(\mu).
\]
Since $\mu$ has a compact support, there exists $R>0$ such that
\[
|m_k(\mu)|\le R^k,\qquad k\ge1.
\]
Hence
\[
\begin{aligned}
|r_k(\mu)|
&\le
\sum_{p=1}^k
\sum_{\substack{i(1),\ldots,i(p)\ge1\\
i(1)+\cdots+i(p)=k}}
|m_{i(1)}(\mu)|\cdots |m_{i(p)}(\mu)| \\
&\le
R^k
\sum_{p=1}^k
\#\Bigl\{
(i(1),\ldots,i(p)):
i(j)\ge1,\,
i(1)+\cdots+i(p)=k
\Bigr\}.
\end{aligned}
\]
The number of compositions of $k$ into $p$ positive parts is
$\binom{k-1}{p-1}$, and therefore
\[
|r_k(\mu)|
\le
R^k
\sum_{p=1}^k
\binom{k-1}{p-1}
=
R^k2^{k-1}\le (2R)^k.
\]
Since $K_{\mu_t}(z)=
\sqrt t\, K_\mu(\sqrt t\, z),$ comparison of coefficients in \eqref{eq:K-transform generating function} yields
\[
r_k(\mu_t)
=
t^{\frac{2-k}{2}}\,r_k(\mu),
\qquad k\ge1.
\]
Using again the moment-cumulant formula, we obtain
\[
\begin{aligned}
m_{2k}(\mu_t)
&\le
\sum_{p=1}^{2k}
t^{p-k}
\sum_{\substack{i(1),\ldots,i(p)\ge1\\
i(1)+\cdots+i(p)=2k}}
|r_{i(1)}(\mu)|\cdots |r_{i(p)}(\mu)| \\
&\le
(2R)^{2k}
\sum_{p=1}^{2k}
t^{p-k}
\binom{2k-1}{p-1}=(2R)^{2k}t^{1-k}(1+t)^{2k-1}
\end{aligned}
\]
Therefore
\[
\sup_{k\ge1}|m_{2k}(\mu_t)|^{\frac{1}{2k}}
\le
2R\frac{1+t}{\sqrt{t}},
\]
which implies that $\mu_t$ is compactly supported. We complete the proof.
\end{proof}

Next, we provide a proof for Lemma \ref{lem:tech lemma R_delta,varepsilon}.
\begin{proof}
    Recall that $f_\delta(x)=\log |x+i\delta|$ and $g_\delta(x)=xf_\delta(x)$, for a compactly supported probability measure $\mu$ and $z=x+i\varepsilon$ with $\varepsilon>0$,
    \begin{align*}
        R_{\delta,\varepsilon}(\mu) &= \int \left[2\pi u_1(z)v_1(z) +\frac{1}{2} v_1'(z)\right]g_\delta(x)\,dx - \int \frac{1}{2} v_1(z)f_\delta(x)\,dx \\
&\quad + \varepsilon \int \left[ \pi(u_1(z)^2-v_1(z)^2)+\frac{1}{2}u_1'(z)\right]f_\delta(x)\,dx.
    \end{align*}
    As $\varepsilon\to 0$, we aim to show that
    \begin{align}\label{eq:R_delta,varepsilon varepsilon to 0}
        R_{\delta,\varepsilon}(\mu) = \iint &\frac{P_\varepsilon(g_\delta)(y) - P_\varepsilon(g_\delta)(y')}{y-y'}\,d\mu(y)d\mu(y')\nonumber\\
        &- \int \left(f_\delta(x) +\frac{1}{2}\frac{x^2}{x^2+\delta^2}\right)v(x+i\varepsilon)\,dx+O_\delta(\varepsilon).
    \end{align}
    To be more specific, we shall analyze each term of $R_{\delta,\varepsilon}(\mu)$ separately, we still denote $\pi u(z)=\Re G_\mu(z),$ and $\pi v(z)=-\Im G_\mu(z)$ with $z=x+i\varepsilon$.

\medskip
\noindent\textbf{The derivative term.}
An integration by parts gives
\begin{align*}
    \frac{1}{2}\int v_1'(z)\cdot g_\delta(x)\,dx
    &=-\frac{1}{2}\int v_1(z)\cdot g_\delta'(x)\,dx\\
    &=-\frac{1}{2}\int \left(f_\delta(x)+\frac{x^2}{x^2+\delta^2}\right)v_1(z)\, dx.
\end{align*}
The boundary term vanishes because $|v(z)|=O(|z|^{-2})$ for large $|z|$. This computation, combined with the term $-\frac{1}{2}\int f_\delta(x)v(z)\,dx $ gives
\begin{equation}\label{eq:derivative v terms}
    -\int \left(f_\delta(x)+\frac{1}{2}\frac{x^2}{x^2+\delta^2}\right)v(x+i\varepsilon)\,dx.
\end{equation}
Similarly, since $|u(z)|=O(|z|^{-1})$ for large $|z|$, we have
\begin{align*}
    \int f_\delta(x)\cdot u_1'(z)\,dx 
    &= -\int f_\delta'(x)\cdot u_1(z)\,dx\\
    &=-\frac{1}{\pi}\iint f_\delta'(x)\frac{x-y}{(x-y)^2+\varepsilon^2}\,dx\,d\mu(y)\\
    &=\int H_\varepsilon(f_\delta')(y)\,d\mu(y),
\end{align*}
where $H_\varepsilon(x)=\frac{1}{\pi}\frac{x}{x^2+\varepsilon}$ represents the the Hilbert kernel and $H_\varepsilon(f)=H_\varepsilon\ast f$.  Using the residue theorem, one computes
\begin{align}\label{eq:residue}
    H_\varepsilon(f_\delta')(y) &= \Re\{2\pi i[\mathrm{Res}(f_\delta'(z)H_\varepsilon(z-y),i\delta)+\mathrm{Res}(f_\delta'(z)H_\varepsilon(z-y),y+i\varepsilon)]\}\nonumber\\
    &= -\frac{\delta+\varepsilon}{y^2+(\delta+\varepsilon)^2}.
\end{align}
Hence,
\begin{equation}\label{eq:derivative u term upper bound}
    \left|\int f_\delta(x)\cdot u_1'(z)\,dx\right|
= \left|\int \frac{\delta+\varepsilon}{y^2+(\delta+\varepsilon)^2}\,d\mu(y)\right|
\le \frac{1}{\delta},
\end{equation}

\medskip
\noindent\textbf{The quadratic term.} Recall that $P_\varepsilon(x)=\frac{1}{\pi}\frac{\varepsilon}{x^2+\varepsilon^2}$ is the Poisson kernel and $P_\varepsilon(g)=P_\varepsilon\ast g$. By Fubini's theorem,
\begin{align}\label{eq:quadratic uv term}
    \int 2\pi u(z)v(z)\cdot g_\delta(x)\,dx
    &=-\frac{1}{\pi}\Im \int g_\delta(x)\cdot G_\mu(z)^2\,dx\nonumber\\
    &= -\frac{1}{\pi}\Im\iint \frac{1}{y-y'}\int \left[ \frac{g_\delta(x)}{z-y} - \frac{g_\delta(x)}{z-y'} \right]dx\,d\mu(y)d\mu(y')\nonumber\\
    &= \iint \frac{P_\varepsilon( g_\delta)(y)- P_\varepsilon(g_\delta)(y')}{y-y'}\,d\mu(y)d\mu(y').
\end{align}
Since $\pi (u(z)^2-v(z)^2)=\Re (G_\mu(z)^2),$ a similar computation yields
$$
\int f_\delta(x)\cdot \pi(u(z)^2-v(z)^2)\,dx = \iint \frac{H_\varepsilon( f_\delta)(y)- H_\varepsilon(f_\delta)(y')}{y-y'}\,d\mu(y)d\mu(y').
$$
Using \eqref{eq:residue}, observe that
\begin{align*}
    H_\varepsilon(f_\delta)'(x) &= \frac{1}{\pi}\int f_\delta(y)\frac{\partial}{\partial x}\left(\frac{x-y}{(x-y)^2+\varepsilon^2}\right)\,dy\\
    &=-\frac{1}{\pi}\int f_\delta(y)\frac{\partial}{\partial y}\!\left(\frac{x-y}{(x-y)^2+\varepsilon^2}\right)\,dy\\
    &= \frac{1}{\pi}\int f_\delta'(y)\cdot \frac{x-y}{(x-y)^2+\varepsilon^2}\,dy\\
    &= H_\varepsilon(f_\delta')(x) = - \frac{\delta+\varepsilon}{x^2+(\delta+\varepsilon)^2},
\end{align*}
which implies that $\|H_\varepsilon(f_\delta)'\|_{\infty}\le \frac{1}{\delta}$, and thus $H_\varepsilon(f_\delta)$ is $\frac{1}{\delta}$-Lipschitz for any $\varepsilon>0$, we deduce that
\begin{align}\label{eq:quadratic u^2-v^2 term upper bound}
    \left|\int f_\delta(x)\cdot \pi(u(z)^2-v(z)^2)\,dx\right| \le  \frac{1}{\delta}.
\end{align}

Combining \eqref{eq:derivative v terms},\eqref{eq:derivative u term upper bound},\eqref{eq:quadratic uv term} and \eqref{eq:quadratic u^2-v^2 term upper bound}, we obtain the desired equality \eqref{eq:R_delta,varepsilon varepsilon to 0}. Now let $\varepsilon\to 0$, due to the weak convergence of $v(\cdot+i\varepsilon)\to \mu$ on any compact support as $\varepsilon\to 0$ (see e.g. \cite[Theorem 2.4.3]{anderson2010introduction}), we deduce that
\begin{equation}\label{eq:entropy derivative term varepsilon to 0}
    \lim_{\varepsilon\to 0}\int \left(f_\delta(x)+\frac{1}{2}\frac{x^2}{x^2+\delta^2}\right)v(x+i\varepsilon)\,dx = \int \left(f_\delta(x)+\frac{1}{2}\frac{x^2}{x^2+\delta^2}\right)\,d\mu(x).
\end{equation}
To deal with the quadratic term of $R_{\delta,\varepsilon}(\mu_s)$, note that $P_\varepsilon( g_\delta)(x) \to g_\delta(x)$ as $\varepsilon\to 0$, uniformly on the compact support of $\mu$. Furthermore, since $|g_\delta'(x)|$ is bounded by some constant $C_\delta$ on the compact support of $\mu$, we obtain the following bound:
\begin{align*}
    \left|P_\varepsilon(g_\delta)'(x)\right|
    &= \left|\frac{1}{\pi}\int g_\delta(y)\frac{\partial}{\partial x}\!\left(\frac{\varepsilon}{(x-y)^2+\varepsilon^2}\right)\,dy\right|\\
    &= \left|\frac{1}{\pi}\int g_\delta(y)\frac{\partial}{\partial y}\!\left(\frac{\varepsilon}{(x-y)^2+\varepsilon^2}\right)\,dy\right|\\
    &= \left|\frac{1}{\pi}\int g_\delta'(y)\frac{\varepsilon}{(x-y)^2+\varepsilon^2}\,dy\right|\le C_\delta.
\end{align*}
Thus, for any $y,y'\in \mathrm{supp}(\mu)$
\[
\left|\frac{P_\varepsilon(g_\delta)(y) - P_\varepsilon(g_\delta)(y')}{y-y'}\right| \le C_\delta,
\]
and the dominated convergence theorem, together with \eqref{eq:entropy derivative term varepsilon to 0} imply that, as $\varepsilon\to 0$,
\[
R_{\delta,\varepsilon}(\mu)\longrightarrow
\iint \left(\frac{g_\delta(x)-g_\delta(y)}{x-y}
-\frac{1}{2}f_\delta(x)-\frac{1}{2}f_\delta(y)
-\frac{1}{2}\frac{x^2}{x^2+\delta^2} \right)\,d\mu(x)d\mu(y).
\]

This completes the proof of the first part as $\varepsilon\to 0$. To show the monotone convergence of $R_\delta(\mu)$ \eqref{eq:limit R_delta delta to 0}, denote
\[
F_\delta(x,y) = \frac{g_\delta(x)-g_\delta(y)}{x-y}
-\frac{1}{2}f_\delta(x)-\frac{1}{2}f_\delta(y)
\]
and obviously, as $\delta\to 0$, we have $F_\delta(x,y)\to F_0(x,y)$ for all $x,y\in\mathrm{supp}(\mu)$. Moreover, note that
\begin{align*}
    \frac{\partial}{\partial\delta}F_\delta(x,y) &=\frac{\delta}{x-y}\left(\frac{x}{x^2+\delta^2}-\frac{y}{y^2+\delta^2}\right)-\frac{1}{2}\left(\frac{\delta}{x^2+\delta^2}+\frac{\delta}{y^2+\delta^2}\right)\\
    &=-\frac{\delta}{2}\frac{(x+y)^2}{(x^2+\delta^2)(y^2+\delta^2)}\le 0,
\end{align*}
which implies that $F_\delta(x,y)\nearrow F_0(x,y)$ as $\delta\to 0$. It remains to prove that $R_\delta\ge 0$, notice that
\[
F_\delta(x,y)
=
\frac{x+y}{x-y}
\log\left|\frac{x+i\delta}{y+i\delta}\right|
\]
and
\[
\log\left|\frac{x+i\delta}{y+i\delta}\right|
=
\frac12\log\left(\frac{x^2+\delta^2}{y^2+\delta^2}\right)
\]
has the same sign as $x^2-y^2=(x+y)(x-y)$, it follows that $F_\delta(x,y)\geq 0$ for all $x\neq y$. On the diagonal $\{x=y\}$, the singularity is removable, and by the continuity we deduce that $F_\delta\geq 0$ on $\mathbb{R}^2$. Finally, the monotone convergence theorem implies that
\[
\lim_{\delta\downarrow 0}\iint F_\delta(x,y)\,d\mu(x)d\mu(y)
= \iint F_0(x,y)\,d\mu(x)d\mu(y),
\]
which completes the proof.
\end{proof}
In summary, up to a scaling, we first extend the discrete Boolean convolution powers $\mu^{\uplus n}$ to $\mu^{\uplus t}$ for all real $t\ge 1$, and then establish the monotonicity of $\Gamma$ along the continuous family $(\mu^{\uplus t})_{t\ge 1}$. This is in analogy with the monotonicity of free entropy along the free convolution powers $\{\mu^{\boxplus t}\}_{t\ge 1}$, proved by Shlyakhtenko and Tao \cite{shlyakhtenko2020fractional}.

Moreover, we have shown that the map
\[
t\longmapsto \Gamma(\mu_t)
=
\int \log|x|\,d\mu_t(x)
\]
is differentiable whenever \(\mu\) is compactly supported and
\(\Gamma(\mu)>-\infty\). In fact, such a result remains valid for a larger class of functions under suitable
regularity assumptions.

\begin{corollary}\label{cor:general h along boolean clt}
    Let \(h\in C^1(\mathbb{R})\) be such that
    \(|h(x)|=o(|x|)\) as \(|x|\to\infty\), and \(H_\varepsilon(h')\) is bounded on every compact set uniformly in \(\varepsilon\). Define
    \[
    T_h(\mu)=\int h(x)\,d\mu(x)
    \]
    for any compactly supported probability measure $\mu$ and
    \[
    D_h(\mu)=\lim_{\Delta t\to0}
    \frac{T_h(\mu_{1+\Delta t})-T_h(\mu)}{\Delta t}.
    \]
    Then, for any $t\ge 1$, we have
    \[
    T_h(\mu_t) - T_h(\mu)
    =\int_1^s \frac{1}{s}D_h(\mu_{t_0})\,ds.
    \]
    Moreover, \(D_h(\mu)\) admits the explicit expression
    \begin{equation*}
    D_h(\mu)
    =
    \iint
    \frac{xh(x)-yh(y)}{x-y}\,d\mu(x)d\mu(y)
    -\int h(x)\,d\mu(x)
    -\frac12\int xh'(x)\,d\mu(x).
    \end{equation*}
    In particular, the function \(t\mapsto T_h(\mu_t)\) is
    $C^1$.
\end{corollary}

\begin{remark}\label{rem:existence of other functional for monotonicity}
    The function \(\log|\cdot|\) does not belong to the class $C^1$ due to its singularity at the origin. Therefore, an additional argument is required to justify the passage to the limit from \(\log|x+i\varepsilon|\in C^1(\mathbb{R})\), as in Lemma \ref{lem:tech lemma R_delta,varepsilon}. As a consequence of the above corollary, if there exists a function \(h\) such that $D_h(\mu)\ge 0$ for any $\mu\in\mathcal{P}_{0,c}^2$,
    then
    \[
    T_h(\mu_n)=\int h(x)\,d\mu_n(x)
    \]
    is also non-decreasing along the Boolean CLT. Indeed, one can check that the function $f_\delta(x)=\log|x+i\delta|$ satisfies this condition. Hence, the monotonicity property does not characterize the proposed Boolean entropy uniquely; other functionals may enjoy the same monotonicity behavior.
\end{remark}

We also point out that our arguments for the monotonicity \eqref{ineq:mono entropy} extend to any compactly supported probability measures, without assuming centeredness. More generally, for any $\mu\in\mathcal{P}^2$, one can apply a truncation argument to extend the monotonicity properties.
\begin{corollary}
    If $\mu\in \mathcal{P}^2$ and $\Gamma(\mu)>-\infty$, then we have
    $$
    \gamma_\mu(n)\ge \gamma_\mu(m)
    $$
    for any $n\ge m\ge 1$.
\end{corollary}
\begin{proof}
    For any $\mu\in\mathcal{P}^2$, the $K$-transform admits the expansion
\[
K_{\mu_t}(z)
=
\sqrt t\,K_\mu(\sqrt t\,z)
=
\sqrt t\,r_1(\mu)
+\frac{r_2(\mu)}{z}
+o(z^{-1}),
\qquad z\to\infty.
\]
    Thus, $\mu_t\in\mathcal{P}^2$ for every $t\ge 1$ and it is enough to prove that
\[
\Gamma(\mu)\le \Gamma(\mu_t),
\qquad
\mu\in\mathcal P^2,\ t\ge1.
\]
Indeed, by the multiplicative property,
\[
\Gamma(\mu_t)
=
\Gamma\bigl((\mu_s)_{t/s}\bigr)
\ge
\Gamma(\mu_s),
\qquad
t\ge s\ge1.
\]

To this end, for \(R>0\), define the truncated probability measure
\[
d\mu_R(x)
=
\mu([-R,R])^{-1}\1_{[-R,R]}(x)\,d\mu(x).
\]
In general, \(\mu_R\) is not centered.
Nevertheless, the monotonicity \eqref{ineq:mono entropy} remains valid, namely,
\[
\Gamma((\mu_R)_t)
\ge
\Gamma(\mu_R),
\qquad
R>0,\ t\ge1.
\]
Since weak convergence of probability measures is equivalent to the pointwise
convergence of their Cauchy transforms on $\mathbb{C}_+$,
\(\mu_R\to\mu\) implies that
\(G_{\mu_R}\to G_\mu\) pointwise. It then follows from
\eqref{eq:boolean CLT cauchy flow} that $(\mu_R)_t\to \mu_t.$ Consequently, $K_{(\mu_R)_t}(z)
\to
K_{\mu_t}(z),$
so that $r_i((\mu_R)_t)\to r_i(\mu_t),\ i=1,2,$ which yields
\[
m_2((\mu_R)_t)\longrightarrow m_2(\mu_t).
\]

On the other hand, the functional
\[
\nu\longmapsto
\int\bigl(x^2-\log|x|\bigr)\,d\nu(x)
\]
is lower semicontinuous with respect to the weak topology. Therefore,
\[
\liminf_{R\to\infty}
\int
(x^2-\log|x|)
\,d(\mu_R)_t(x)
\ge
\int
(x^2-\log|x|)
\,d\mu_t(x).
\]
Combining this with the convergence of second moments, we obtain
\[
\limsup_{R\to\infty}
\Gamma((\mu_R)_t)
\le
\Gamma(\mu_t).
\]
Finally, by the dominated convergence theorem and \(\mu([-R,R])\to1\),
\[
\Gamma(\mu_R)
\longrightarrow
\Gamma(\mu).
\]
Therefore,
\[
\Gamma(\mu)
=
\lim_{R\to\infty}
\Gamma(\mu_R)
\le
\limsup_{R\to\infty}
\Gamma((\mu_R)_t)
\le
\Gamma(\mu_t),
\]
which completes the proof.
\end{proof}
The monotonicity in Theorem \ref{thm:max and mono} then follows from the above Corollary.

\appendix
\renewcommand{\thetheorem}{\Alph{section}.\arabic{theorem}}
\section{Large deviation principle}\label{LDP appendix}
This appendix recalls basic definitions and main results of the large deviation theory. We refer to the readers to \cite{deuschel2001large} and \cite{dembo2009large} for a full treatment.

In what follows, $X$ will be assumed to be a Polish space (that is a complete separable metric space). 

\begin{definition}\label{weak LDP}
    A sequence $(\mu_N)_{N\in\mathbb{N}}$ of probability measure on $X$ satisfies a weak large deviation principle if (a) and (b) hold, and in addition (c) holds for all compact sets $F\subset X.$
\end{definition}

The proof of the large deviation principle often proceeds first by the proof of a weak large deviation principle, in conjunction with the so-called exponential tightness property.

\begin{definition}\label{def exponential tightness}
    A sequence $(\mu_N)_{N\in\mathbb{N}}$ of probability measure on $X$ is exponentially tight iff there exists a sequence $(K_L)_{L\in\mathbb{N}}$ of compact sets such that
\begin{equation}
    \limsup_{L\to \infty}\limsup_{N\to \infty}\frac{1}{a_N}\log \mu_N(K_L^c)=-\infty.
\end{equation}
\end{definition}

The interests in these concepts lie in the following
\begin{theorem}\label{theorem LDP}
If $(\mu_N)_{N\in \N}$ satisfies the exponential tightness and weak large deviation principle, then $(\mu_N)_{N\in\N}$ satisfies the large deviation principle.
\end{theorem}

The following lemma provides a simpler way to prove the weak large deviation principle.

\begin{lemma}\label{weak LDP lemma}
Assume that for all $x\in X,$ 
\begin{equation}
    \lim_{\varepsilon\to 0}\limsup_{N\to \infty}\frac{1}{a_N}\log\mu_N(B(x , \varepsilon))=\lim_{\varepsilon\to 0}\liminf_{N\to\infty}\frac{1}{a_N}\log\mu_N(B(x , \varepsilon))=-I(x).
\end{equation}
then $(\mu_N)_{N\in\N}$ satisfies the weak large deviation principle.
\end{lemma}

\bibliographystyle{plain}
\bibliography{reference}

\end{document}